\documentclass{article}

\usepackage{graphicx}
\usepackage{indentfirst}
\usepackage{booktabs}
\usepackage[ruled,linesnumbered]{algorithm2e}
\usepackage{lineno}
\usepackage{xcolor,hyperref}
\usepackage{tikz}
\usepackage[utf8]{inputenc}
\usepackage{amsfonts}
\usepackage{multirow}
\usepackage{lipsum}
\usetikzlibrary{matrix, positioning, calc}
\usepackage{csquotes}
\usepackage{amsfonts,amsmath,amssymb}

\usepackage[margin=1in]{geometry} 
\usetikzlibrary{arrows, arrows.meta}
\usepackage{authblk}
\newcommand{\keywords}[1]{\textbf{\textit{Keywords:}} #1}

\begin{document}

\title{Towards the efficient calculation of quantity of interest from
  steady Euler equations II: a CNNs-based automatic implementation.}

\author[1]{Jingfeng Wang}
\author[1,2,3]{Guanghui Hu\thanks{Corresponding author: \texttt{garyhu@um.edu.mo}}}
\affil[1]{Department of Mathematics, Faculty of Science and Technology, University of Macau, Macao S.A.R., China}
\affil[2]{Zhuhai UM Science and Technology Research Institute, Zhuhai, Guangdong, China}
\affil[3]{Guangdong-Hong Kong-Macao Joint Laboratory for Data-Driven Fluid Mechanics and Engineering Applications, University of Macau, China}

\maketitle

\begin{abstract}
  In \cite{wang2023towards}, a dual-consistent
  dual-weighted residual-based $h$-adaptive method has been proposed
  based on a Newton-GMG framework, towards the accurate calculation of a
  given quantity of interest from Euler equations. The performance of
  such a numerical method is satisfactory, i.e., the stable
  convergence of the quantity of interest can be observed in all
  numerical experiments. In this paper, we will focus on the
  efficiency issue to further develop this method, since efficiency is
  vital for numerical methods in practical applications such as the
  optimal design of the vehicle shape. Three approaches are studied
  for addressing the efficiency issue, i.e., i). using convolutional
  neural networks as a solver for dual equations, ii). designing an
  automatic adjustment strategy for the tolerance in the $h$-adaptive
  process to conduct the local refinement and/or coarsening of mesh
  grids, and iii). introducing OpenMP, a shared memory parallelization
  technique, to accelerate the module such as the solution
  reconstruction in the method. The feasibility of each approach and
  numerical issues are discussed in depth, and significant
  acceleration from those approaches in simulations can be observed
  clearly from a number of numerical experiments. In convolutional
  neural networks, it is worth mentioning that the dual consistency
  plays an important role to guarantee the efficiency of the whole
  method and that unstructured meshes are employed in all
  simulations.
\end{abstract} 

\keywords{DWR-based $h$-adaptivity; Dual consistency; Convolutional neural networks; Automatic adjustment of tolerance; Newton-GMG method for Euler equations}

\maketitle


\section{Introduction}

In \cite{wang2023towards}, a dual-consistent dual-weighted
residual(DWR)-based adaptive mesh method has been constructed, from
which the smooth convergence of quantity of interest can be
observed. This technique is particularly advantageous in applications
such as the optimal design of vehicle shapes where the quantity of
interest, e.g., the lift-to-drag ratio, needs an accurate evaluation.

It is noted that the efficiency of a numerical method is a key feature
towards the practical applications, as in a typical problem of optimal
design of the vehicle shape, governing partial differential equations
(PDEs) need to be solved numerous times, depending on factors such
as the number of design parameters, mesh resolution, optimization
method\cite{jameson2003aerodynamic}. As a consequence, in this paper, efforts are devoted to
further improving the efficiency of the numerical method proposed in
\cite{wang2023towards}, to make the method and related library
AFVM4CFD a competitive one in the market.

For such a purpose, three approaches will be studied in depth in this
paper, including i). using convolutional neural networks (CNNs) to
produce numerical solutions of the dual problem, ii). designing an
automatic adjustment strategy for the tolerance in the $h$-adaptive
process, and iii). introducing OpenMP to parallelize modules such as
the solution reconstruction and dual solver in the method. In terms of where issues
are from, the first approach is for delivering more efficient
solutions of dual equations, while the latter two are for accelerating
the algorithm.

It is noted that Galerkin orthogonality between numerical solutions
and related error is the reason why Hartmann\cite{hartmann2007adjoint}
suggested solving dual equations in a larger space. In
\cite{wang2023towards}, the refinement of mesh grids is employed for
such a purpose. However, this approach is not a good choice from the
efficiency point of view\cite{liu2022mp}. Although there have been
other choices such as raising the order of approximate
polynomials\cite{dolejsi2022}, high order interpolation of numerical solutions of
dual equations\cite{venditti2000adjoint}, it should be pointed out that, i). solving dual
equations using a classical solver involved in all aforementioned
approaches, in which the time-consuming modules such as the construction
of the finite-dimensional space, solving the system of linear
equations, etc., need to be implemented, and ii). numerical solutions
of dual equations only serve as the weight in the DWR-based error
estimation, which means that the numerical accuracy in this scenario
is not the first concern.

CNNs, based on the above two observations, become an ideal choice for the
solver of dual equations, due to its ability on fastly generating
acceptable numerical solutions of dual equations. In the recent past,
neural networks have been successfully applied for implementing the
DWR-based h-adaptation method in computational fluid
dynamics(CFD). For instance, in \cite{chen2021output, chen2020output},
an encoder-decoder algorithm is developed to predict the indicators
generated from dual equations, and in \cite{chakraborty2021multigoal},
primal solutions and dual solutions are both obtained from deep neural
networks. Similarly, a data-driven goal-oriented mesh adaptation
approach is developed to generate the error indicators in
\cite{wallwork2022e2n}. The recent breakthrough in convolutional
neural networks make training in high-dimensional problems
possible to handle numerical computations in CFD. Besides, libraries
build in Tensorflow\cite{abadi2016tensorflow},
Pytorch\cite{paszke2017automatic}, and cuDNN provide powerful tools to
design specific models that cater to the demands of CFD issues, where
the GPUs module can be adopted for further improvement in the training
process. Another reason to choose CNNs in this work is its
generalization ability, i.e., quality numerical solutions of dual
equations with different attack angles, Mach numbers, etc. can be
generated effectively. Furthermore, CNNs are not sensitive to the
geometry of the domain, which is also a feature desired.

To accelerate the numerical simulation, the
first approach in this paper is to design an automatic adjustment
strategy for the tolerance used in the $h$-adaptive process. It is
known that in an $h$-adaptive mesh method, after an error indicator is
generated for each element, a tolerance needs to be given to guide the
local refinement and/or coarsening of mesh grids so that the global
error of numerical solutions can be effectively controlled. In
\cite{wang2023towards}, the strategy for setting up this tolerance is
to use a fixed and sufficiently small one during the whole
simulation. This is not a good strategy since in the first several
rounds of adaptive refinement of mesh grids, unreliable error
indicators would result in over-refinement and/or coarsening of mesh
grids if a sufficiently small tolerance is employed. Furthermore, the
nonlinearity of governing equations also implies that a dynamic
adjustment of the tolerance would be a better choice. Motivated by the
decreasing threshold method proposed in
\cite{nemec2014toward,nemec2007adjoint}, error indicators are plotted
for analysis and the corresponding dynamic tolerance selection method
is developed in this work. Then an automatic DWR-based h-adaptivity
method is developed without manual intervention.

The second approach used in this paper for accelerating simulations is
to use OpenMP \cite{chandra2001parallel} to parallelize modules both for solving
primal and dual equations. It is noted that the hardware used in this work for
numerical experiments is a workstation with multiple CPU cores, which
makes OpenMP a perfect choice for parallelization. In OpenMP, data
race is a key issue that should be avoided during the parallelization
of the algorithm. The good news is that in our Newton-GMG framework,
there are many time-consuming modules, such as the solution
reconstruction, the update of cache information for each element, as well
as the multigrid method, which can be implemented without data race. Hence,
in this paper, the OpenMP parallelization of the algorithm and code
proposed in \cite{wang2023towards} will be discussed in detail for the
acceleration.

In this paper, we will follow the aforementioned three approaches to
further improve the efficiency of the numerical method proposed in
\cite{wang2023towards}. Firstly, we have constructed a CNNs-based dual solver. For fear of overfitting, techniques such as batch normalization and pooling have been incorporated. Furthermore, in order to circumvent the issues of gradient vanishing and explosion, initialization and residual connection modules are integrated, as highlighted in references \cite{he2015delving,he2016deep}. This enhanced solver efficiently predicts the dual equations on the current mesh, facilitating rapid adaptation processes. Notably, the trained model demonstrates generalization capabilities across training datasets, configurations, and multiple airfoils. Secondly, an automatic adjustment strategy for the tolerance is devised in the $h$-adaptive process. This strategy revolves around plotting the distribution of mesh adaptation steps and analyzing them using the Kolmogorov-Smirnov method. Our findings suggest that the distribution of these indicators closely follows the Weibull distribution, characterized by a KS-stat at its uppermost value. However, it should be noted that the parameters might differ for various adaptation steps. In some configurations, the distribution even aligns more closely with the gamma distribution. Owing to these variations, relying on a mathematical expression to select tolerance values is deemed impractical. Instead, we have employed a methodology where values of indicators are systematically sorted, and a fixed proportion is chosen to derive the tolerance for different adaptation steps. Finally, three modules
  in our algorithm, i.e., solution reconstruction, update cache
  information for each element and CNNs form dual solver are parallelized with OpenMP. The scalability of these modules is illustrated in experiments with the different number of threads. Using the standard metric of speedup, typical in parallel computing evaluations, our results indicate a substantial enhancement in the performance of the framework.

It is found from the numerical experiments that the CNNs dual solver
saves time in an order of magnitude compared with the traditional
solver. More specifically, with time complexity of $O(n)$. Besides,
the dynamic tolerance strategy accelerates the calculation of the
target functional and benefits for automatic selection without manual
intervention. With such a strategy, a stable growth rate of mesh size
can be expected. In \cite{wang2023towards}, the importance of dual
consistency within the DWR-based h-adaptation is discussed, such
property extends to the construction of training datasets, which has
been further elaborated upon in this research. The trained model has
generalization capabilities that configurations beyond the trained
category can still be simulated precisely. It is worth noting that our
CNNs solver can generate reliable dual solutions on unstructured
mesh. With the introduced parallel module, the time spent on primal solver and dual solver has been saved a lot.

The rest of this paper is organized as follows. In Section 2, we
present the fundamental notations and provide a concise introduction
to the Newton-GMG solver. In Section 3, the architecture of our
convolutional neural networks is introduced and the training
performance is discussed. We emphasized the importance of dual
consistency in the training process. In section 4, two acceleration strategies are discussed. We organize the
dynamic tolerance selection method and analyze the effects of this
algorithm. Then parallel computing is adopted. In section 5, we present the details of the numerical
experiments, demonstrating the reliability of our CNNs-based dual
solver.

\section{A Brief Introduction to Dual-Consistent DWR-based h-adaptive Newton-GMG}
\subsection{Basic notations}
Let $\Omega$ be a domain in
$\mathbb R^2$ with boundary $\Gamma$. The inviscid two-dimensional steady Euler equations derived from the conservation law can be written as:
Find $\mathbf{u}~:~\Omega\rightarrow \mathbb R^4 $ such that~
\begin{equation}\label{primal}
  \nabla\cdot \mathcal{F}( \mathbf{u})=0, \quad \mbox{in}~~\Omega,
  \end{equation}
subject to a certain set of boundary conditions. 
Here $\mathbf{u}$ and $\mathcal F(\mathbf{u})$ are the conservative variables and flux, which are given by \begin{equation}
    \mathbf{u}=\begin{bmatrix}
      \rho \\ \rho u_x \\\rho u_y\\E
    \end{bmatrix},
    \qquad\text{and}~\mathcal F(\mathbf{u})=\begin{bmatrix}
      \rho u_x & \rho u_y
      \\ \rho u_x^2+p&\rho u_xu_y
      \\ \rho u_xu_y & \rho u_y^2+p
      \\ u_x(E+p) & u_y(E+p)
    \end{bmatrix},
\end{equation}  
where $(u_x,u_y)^T, \rho, p, E$ denote the velocity, density, pressure, and total energy, respectively. We use the equation of state to close the system, which is 
 \begin{equation}
   E=\frac{p}{\gamma-1}+\frac{1}{2}\rho(u_x^2+u_y^2).
 \end{equation}
Here $\gamma=1.4$ is the ratio of the specific heat of the perfect gas.

\subsection{Finite volume discretization}
In \cite{HU2016235}, equations \eqref{primal} are solved with a high-order finite volume method. The methodology will be briefly introduced as follows. To derive a discretization scheme, a shape
regular subdivision $\mathcal{T}$ turns $\Omega$ into different control volumes, $K$. $K_i$ is used to define the $i$-th element in this subdivision. $e_{i,j}$ denotes the common edge of $K_i$ and $K_j$, i.e., $e_{i,j}=\partial K_{i} \cap \partial K_{j}$. The unit outer normal vector on the edge $e_{i,j}$ with respect to $K_{i}$ is represented as $n_{i,j}$. Consequently, with the divergence theorem, the Euler equations in this discretized domain can be reformulated as 
\begin{equation}
    \label{fvm_euler}
    \mathcal{A}(\mathbf{u})=\int_\Omega \nabla \cdot \mathcal{F}(\mathbf{u})dx=\sum\limits_{i}\int_{K_i}\nabla\cdot
    \mathcal{F}(\mathbf{u})dx=\sum\limits_{i}\sum\limits_{j}\oint_{e_{i,j}\in\partial K_i}\mathcal{F}(\mathbf{u})\cdot n_{i,j}ds=0.
\end{equation}
With the numerical flux $\mathcal{H}(\cdot,\cdot,\cdot)$ introduced in this equation, a fully discretized system can be obtained as 
\begin{equation}
    \label{EulerDiscrete}
    \sum\limits_{i}\sum\limits_{j}\oint_{e_{i,j}\in\partial K_i}\mathcal{H}(\mathbf{u}_i,\mathbf{u}_j, n_{i,j})ds=0.
\end{equation}
\subsection{Newton-GMG solver}
To solve the Equation \eqref{EulerDiscrete}, the Newton method is employed for the linearization and a linear multigrid method proposed in \cite{li2008multigrid} is utilized for solving the system. At first, the Equation \eqref{EulerDiscrete} is expanded by the Taylor series. Then we neglect the higher-order terms. Then the system becomes

\begin{equation}
    \begin{aligned}
        \label{regularized_equation}
        \displaystyle \alpha \left|\!\left|\sum\limits_{i}\sum\limits_{j}\int_{e_{i,j}\in\partial\mathcal{K}_i}\mathcal{H}(\mathbf{u}_i^{(n)},\mathbf{u}_j^{(n)}, n_{i,j})ds \right|\!\right|_{L_1}\Delta \mathbf{u}_i^{(n)}&+\sum\limits_{i}\sum\limits_{j}\int_{e_{i,j}\in\partial\mathcal{K}_i}\Delta \mathbf{u}_i^{(n)}\frac{\partial\mathcal{H}(\mathbf{u}_i^{(n)},\mathbf{u}_j^{(n)}, n_{i,j})}{\partial \mathbf{u}_i^{(n)}}ds\\
&+\sum\limits_{i}\sum\limits_{j}\int_{e_{i,j}\in\partial\mathcal{K}_i}\Delta \mathbf{u}_j^{(n)}\frac{\partial\mathcal{H}(\mathbf{u}_i^{(n)},\mathbf{u}_j^{(n)}, n_{i,j})}{\partial \mathbf{u}_j^{(n)}}ds\\
&=-\sum\limits_{i}\sum\limits_{j}\int_{e_{i,j}\in\partial\mathcal{K}_i}\mathcal{H}(\mathbf{u}_i^{(n)},\mathbf{u}_j^{(n)}, n_{i,j})ds,
\end{aligned}
\end{equation}
where ${\partial\mathcal{H}(\cdot,\cdot,\cdot)}/{\partial \mathbf{u}}$ denotes the Jacobian matrix of numerical flux, and $\Delta \mathbf{u}_i$ is the increment of the conservative variables in the $i-$th element. After each Newton iteration, the cell average is updated by $\mathbf{u}_i^{(n+1)}=\mathbf{u}_i^{(n)}+\Delta \mathbf{u}_i^{(n)}$. In the simulation, a regularization term is added to the system. 
This regularization is effective since the local residual can quantify whether the solution is close to a steady state. Based on the observation, the solution is far from a steady state at the initial stage. With the iteration processed, the solution gets updated and approaches the steady state where the regularization term gets close to zero. 

To solve the Equation \eqref{regularized_equation}, the geometrical multigrid technique is applied. In \cite{wang2023towards}, the regularization and geometrical multigrid methods are also applied to the dual equations which perform well to get a robust solver for obtaining dual equations.
\subsection{Dual Weighted Residual Mesh Adaptation}
In practical issues, engineers care about how to calculate the quantity of interest efficiency. For instance, in the field of shape optimal design, lift and drag are of main concern. Then the target functionals are considered as 
\begin{equation}\label{lift_and_drag}
    \mathcal J(\mathbf{u}) =\int_{\Gamma}p_{\Gamma}(\mathbf{ u})\mathbf{n}\cdot\beta,
  \end{equation}
  where $p_{\Gamma}$ is the pressure along the boundary of airfoil and $\beta$ in the above formula is given as 
\begin{equation}\beta=\left\{
  \begin{array}{l}
     (\cos\alpha,\sin\alpha)^T/C_{\infty},\text{ for drag calculation}, \\
     (-\sin\alpha,\cos\alpha)^T/C_{\infty},\text{ for lift calculation}.
  \end{array} 
  \right.
\end{equation}

Indeed, the quantity of interest can be calculated with high precision if the solutions from \eqref{primal} are solved accurately. However, the computational cost increases exponentially with uniform refinement. Then dual-weighted residual method offers an opportunity to calculate the quantity of interest efficiency. In \cite{venditti2000adjoint}, discretized dual equations are solved from the perspective of extrapolation, which is easily implemented on a finite volume scheme. A brief introduction to this method will be given in the following section at first.

Suppose $\mathcal{T}_H$ is the partition of the computational domain into a coarse mesh, and $\mathcal{T}_h$ is the partition into fine mesh correspondingly. The primal equations \eqref{primal} are denoted by a residual form, where $\mathbf {u}_H$ denotes the numerical solutions obtained from the discretization on coarse mesh, i.e.
\begin{equation}
\mathcal{R}_H(\mathbf {u}_H) = 0.    
\end{equation}
 The original motivation of this method is to estimate the target functional $\mathcal{J}(\mathbf{u})$ on a fine mesh $\mathcal{J}_h(\mathbf{u}_h)$ without solving the equations on fine mesh. A multiple-variable Taylor series expansion is applied at first, i.e.,
\begin{equation}
    \label{quantity_vector}
    J_h(\mathbf {u}_h) =J_h(\mathbf{u}_h^H)+\left.\frac{\partial J_h}{\partial \mathbf{u}_h}\right|_{\mathbf{u}_h^H}(\mathbf{u}_h-\mathbf{u}_h^H)+\cdots,
  \end{equation}
  here $\mathbf{u}_h^H$ represents the coarse solution $\mathbf{u}_H$ mapped onto the fine space $\mathcal V_h$ via some prolongation operator $I_h^H$,
\begin{equation}
  \mathbf{u}_h^H= I_h^H\mathbf{u}_H.
\end{equation}
Similarly, the equations on the fine mesh can be expanded through the Taylor expansion, i.e.,
\begin{equation}
    \label{residual_vector}
    \mathcal{R}_h(\mathbf{u}_h) =\mathcal{R}_h(\mathbf{u}_h^H)+\left.\frac{\partial \mathcal{R}_h}{\partial \mathbf{u}_h}\right|_{\mathbf{u}_h^H}(\mathbf{u}_h-\mathbf{u}_h^H)+\cdots.
  \end{equation}
  Symbolically, the Jacobian matrix $\left.\frac{\partial \mathcal{R}_h}{\partial \mathbf{u}_h}\right|_{\mathbf{u}_h^H}$ can be inverted to get the error representation,
  \begin{equation}
    \label{vector_error}
    \mathbf{u}-\mathbf{u}_h^H\thickapprox -(\left.\frac{\partial \mathcal{R}_h}{\partial\mathbf{u}_h}\right|_{\mathbf{u}_h^H})^{-1}\mathcal{R}_h(\mathbf{u}_h^H).
  \end{equation}
  Substituting the expression \eqref{vector_error} into \eqref{quantity_vector}, the quantity of interest can be denoted as 
  \begin{equation}
    \label{quantity_error}
    J_h(\mathbf{u}_h) =J_h(\mathbf{u}_h^H)-(\mathbf{z}_h)^T\mathcal{R}_h(\mathbf{u}_h^H),
  \end{equation} 
  where $\mathbf{\mathbf{z}}_h|_{\mathbf{u}_h^H}\it$ is obtained from \textit{Fully discrete dual equations}:

\begin{equation}\label{discrete_vector}
  \left(\left.\frac{\partial \mathcal R_h}{\partial \mathbf{u}_h}\right|_{\mathbf{u}_h^H}\right)^T\mathbf{z}_h=\left(\left.\frac{\partial \mathcal J_h}{\partial \mathbf{u}_h}\right|_{\mathbf{u}_h^H}\right)^T.
\end{equation}
In our earlier work, $\mathbf{z}_h$ is calculated on the embedded mesh. Various strategies have been proposed in the literature to simplify this process. For instance, \cite{venditti2000adjoint}\cite{nemec2007adjoint} used the error correction method which interpolates the dual solutions to the fine mesh, and in \cite{yamahara2017adaptive}, a smoothing technique is proposed to modify the oscillations due to the interpolation. Even though, the mesh needs to be refined at first to make further calculations. In this work, the dual solutions are obtained from the neural networks and can be used to generate the error indicators directly on the current mesh.
The framework of our previous algorithm is briefly summarized as follows.

\begin{algorithm}[H]
    \SetAlgoLined
    \KwData{Initial $\mathcal{K}_H$, $TOL$}
    \KwResult{$\mathcal{K}_h$}
    Using the Newton-GMG to solve $\mathcal{R}_{H}(\mathbf{u}_H)=0$ with residual tolerance $1.0\times 10^{-3}$\;
    Interpolate solution $\mathbf{u}_H$ from the mesh $\mathcal{K}_{H}$ to $\mathcal{K}_{h}$ to get $\mathbf{u}_h^H$\;
    Record the residual $\mathcal{R}_h\left(\mathbf{u}_h^H\right)$\;
    Solve $\mathcal{R}_{h}(\mathbf{u}_h)=0$ with residual tolerance $1.0\times 10^{-3}$\;
    Using the Newton-GMG to solve the dual equation to get $\mathbf{z}_h$\;
    Calculate the error indicator for each element\;
    \While{$\mathcal{E}_{K_H}>TOL$ for some $K_H$}{
      Adaptively refine the mesh $\mathcal{K}_{H}$ with the process in \cite{HU2016235};}
    \caption{DWR for one-step mesh refinement}
  \end{algorithm}
Even though the dual-consistent DWR-based mesh adaptation method shows great potential for accurately resolving quantities of interest, the generation of high-quality mesh continues to be a computationally demanding process. We observe that machine learning techniques, specifically convolutional neural networks (CNNs), may enable the rapid attainment of dual solutions without compromising precision significantly. The application of CNNs in predicting dual solutions only affects the adaptation process, leaving the underlying physical mechanisms of the primal equations unaltered. Consequently, the derivation of dual solutions appears apt for replacement with CNNs-based methodologies.
  In the next section, we will introduce how to train neural networks to get the dual equations so that we can substitute step 4.

  \section{CNNs-based Solutions of Dual Equations}

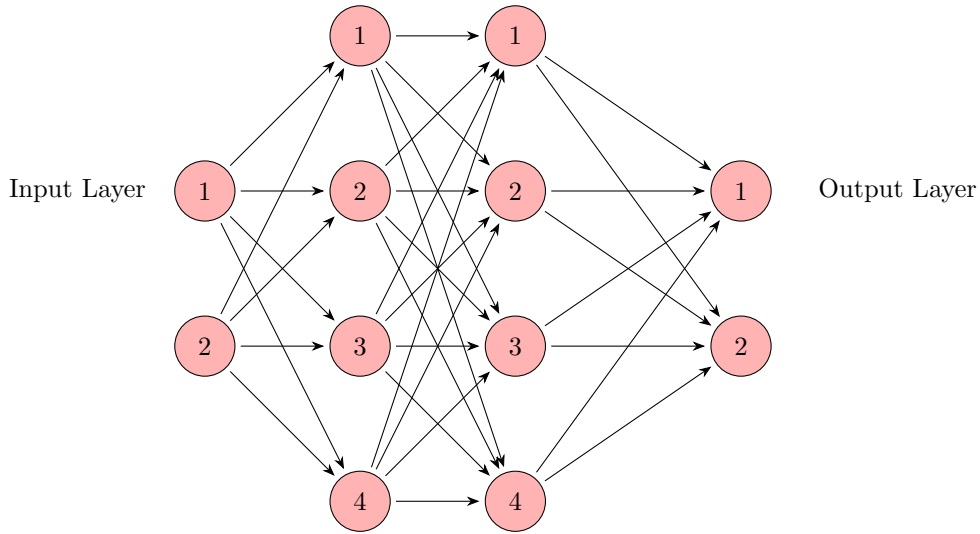
\begin{figure}[!hbtp]
      \centering

  \begin{tikzpicture}[
    node distance = 1.25cm,
    layer/.style={draw, circle, fill=red!30, minimum size=0.8cm},
    arrow/.style={-Stealth, shorten >=2pt, shorten <=2pt}
  ]

  \node[layer] (input1) at (0, -1) {1};
  \node[layer, below=of input1] (input2) {2};
  \node[left=0.25cm of input1] {Input Layer};
  
  \node[layer, right=of input1] (hidden12) {2};
  \node[layer, above=of hidden12] (hidden11) {1};
  \node[layer, below=of hidden12] (hidden13) {3};
  \node[layer, below=of hidden13] (hidden14) {4};

  \node[layer, right=of hidden11] (hidden21) {1};
  \node[layer, below=of hidden21] (hidden22) {2};
  \node[layer, below=of hidden22] (hidden23) {3};
  \node[layer, below=of hidden23] (hidden24) {4};

  \node[layer] (output1) at ([xshift=3cm]hidden22) {1};
  \node[layer, below=of output1] (output2) {2};
  \node[right=0.5cm of output1] {Output Layer};
  
  \foreach \i in {1, 2}
    \foreach \j in {1, 2, 3, 4}
      \draw [arrow] (input\i) -- (hidden1\j);
  
  \foreach \i in {1, 2, 3, 4}
    \foreach \j in {1, 2, 3, 4}
      \draw [arrow] (hidden1\i) -- (hidden2\j);
  
  \foreach \i in {1, 2, 3, 4}
    \foreach \j in {1, 2}
      \draw [arrow] (hidden2\i) -- (output\j);
      \label{fully}
  \end{tikzpicture}
  \caption{Fully connected neural networks.}
  \label{fullyarc}
  \end{figure}

  Convolutional neural networks are very similar to the traditional fully connected neural networks as shown in Figure \ref{fullyarc}, consisting of learnable weights and biases in different layers. However, the practical issues often concern high-dimensional problems where the connection between different neurons and layers may become complicated. In \cite{chen2021output}, autodecoder neural networks are trained to generate indicators and a topology map converts the information to the structured mesh. In this work, since the adaptation is conducted on an unstructured mesh, the map relations are hard to construct, then we build the architecture which learns from the data directly.
  \subsection{Convolutional neural networks}
      In \eqref{discrete_vector}, the dual equations are obtained with the primal solutions. In an abstract form, the dual equations can be treated as a functional of $\mathcal{R}_h, \mathbf{u}_h, \mathcal{J}_h$, i.e,
      \begin{equation}
          \mathbf{z}_h = \mathcal{N}(\mathcal{R}_h, \mathbf{u}_h, \mathcal{J}_h),
      \end{equation}
    where the functional $\mathcal{N}(\cdot,\cdot,\cdot)$ is unknown. Since the information of mesh grids is not included in the training process, $\mathcal{J}_h$ is not suitable to be treated as input. Alternatively, while the quantity of interest is the integration around the boundary of the airfoil in this specific model, the location of the elements and the configurations can work well as inputs. Then the dual solutions can be seen as functional with the expression as follows,
    \begin{equation}
          \mathbf{z}_h = \widetilde{\mathcal{N}}(\mathcal{R}_h, \mathbf{u}_h, \Vec{G}, Ma, \theta),
    \end{equation}
    where $\Vec{G}$ is the coordinates the barycenter of elements, $Ma$ is the Mach number and $\theta$ is attack angle. The neural networks are trained to approximate the unknown functional $\widetilde{\mathcal{N}}(\cdot,\cdot,\cdot,\cdot,\cdot)$. The mechanism behind the relations is complicated and the dimension is high for training a satisfactory model within the fully connected neural networks framework. Then we consider using the CNNs for the training.

  CNNs offer a more efficient approach to handling high-dimensional problems compared to fully connected networks. Instead of connecting all neurons from the preceding layer, CNNs establish connections only within local regions of the input. This specialized architecture consists of multiple learnable filters, each possessing parameters such as weights and biases. When the networks become deeper, the prediction in CNNs is quicker than that in fully connected neural networks. Then we considered building the CNNs model for generating the dual solutions.

  Inputs for CNNs are typically tensors or matrices. As filters slide across different portions of the input, they encode information and transmit it to the subsequent layer. This selective connectivity allows CNNs to exploit spatial hierarchies and local structures in the input data, resulting in more effective and computationally efficient learning. Similarly, a decoder process is also integrated into this framework to pass the encoded information to the final results.
\begin{figure}
    \centering
  \begin{tikzpicture}[mymatrix/.style={matrix of nodes, nodes={draw, minimum size=1cm, anchor=center}, column sep=-\pgflinewidth, row sep=-\pgflinewidth, nodes in empty cells, inner sep=0pt, outer sep=0pt, draw}]
    \matrix[mymatrix] (mat1) {
        1 & 2 & 3 \\
        4 & 5 & 6 \\
        7 & 8 & 9 \\
    };
    
    \matrix[mymatrix, right=of mat1] (mat2) {
        1 & 0 \\
        0 & 1 \\
    };
    
    \matrix[mymatrix, right=of mat2] (mat3) {
        6\\
    };
    \matrix[mymatrix, xshift=1cm, right=of mat3] (mat4) {
        |[text = red]| 7 & 9\\
        13 & 15\\
    };
    
    \node[above=0.5cm of mat1] {Input Matrix};
    \node[above=0.5cm of mat2] {Filter};
    \node[above=0.5cm of mat3] {result};
    \node[above=0.5cm of mat4] {Output Matrix};
    
    \foreach \i in {1}
        \foreach \j in {1}
            \draw[red, thick] ($(mat1-\i-\j.north west)+(-0.1,-0.1)$) rectangle ($(mat1-\the\numexpr\i+1\relax-\the\numexpr\j+1\relax.south east)+(0.1,0.1)$);
    
    \draw[->, thick] (mat1.east) -- (mat2.west);
    \draw[->, thick] (mat2.east) -- (mat3.west);
    \draw[->, thick] (mat3.east) -- (mat4.west);
    \path[->, line width=2pt] (mat3.east) -- (mat4.west) node[midway, above] {b = 1};
\end{tikzpicture}
\caption{Convolutional operation.}
\label{cnnarc}
\end{figure}
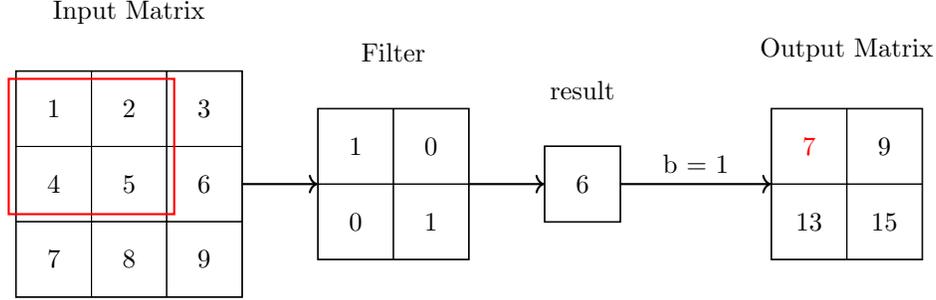

While the training is processed, the parameters will be updated till the preset terminal. In practical issues, connections between inputs and outputs are often hard to capture. By introducing the nonlinear activation function, CNNs can learn more complex patterns and relationships in data. Moreover, the activation function acts as a threshold to determine whether the output could be passed to the next layer. In the backpropagation process, differentiable activation functions, such as the sigmoid, ReLU, or tanh, enable gradient-based learning and help the network learn from the data through iterative weight adjustments. In this work, we adopted the exponential linear unit(ELU) as the activation function\cite{clevert2015fast}, since it has a negative value range for inputs less than zero, which makes the mean activation of the neurons closer to zero. This improvement helps in reducing the bias shif effect significantly compared with other schemes. 

To speed up the training model, we utilized the batch normalization technique, a technique that accelerates convergence by minimizing internal covariate shifts. After this process, a higher learning rate can be set without causing unstable training dynamics. Additionally, for fear of overfitting issues, max pooling and dropout methods are incorporated as well. Max pooling reduces the number of parameters while preserving the most salient features irrespective of their spatial position. Meanwhile, dropout functions as an ensemble technique, simultaneously training multiple subnetworks with diverse architectures. The ultimate output is an average of these subnetworks, resulting in improved generalization. A brief framework of our CNNs can be found in Table \ref{tab:network_architecture}.
\begin{table}[h]
\centering
\begin{tabular}{|l|l|l|l|l|}
\hline
Layer Type & Layer Name & Operation       & Output Channels & Activation \\ \hline
Input      & input         & -               & 3       & -          \\
Conv2d     & enc1       & Encoding        & 64     & -        \\
Inception  & enc2       & Encoding        & 512     & -        \\
Inception  & enc3       & Encoding        & 1024     & -        \\
ResBlock   & enc4       & Encoding        & 512     & ELU        \\
ResBlock   & enc5       & Encoding        & 1024     & ELU        \\
ResBlock   & enc6       & Encoding        & 2048     & ELU        \\
Conv2d     & skip\_adjust & Adjustment        & 2048     & -        \\
ResBlock   & dec1       & Decoding        & 1024      & ELU        \\
Inception  & dec2       & Decoding        & 2048      & -        \\
ResBlock   & dec3       & Decoding        & 1024      & ELU        \\
Inception  & dec4       & Decoding        & 512      & -        \\
AdaptiveAvgPool2d & avg\_pool & Pooling         & 512       & -          \\
Linear     & fc1        & Fully Connected & 64        & ELU        \\
Linear     & fc2\_1        & Fully Connected & 1          & -          \\
Linear     & fc2\_2        & Fully Connected & 1          & -          \\
Linear     & fc2\_3        & Fully Connected & 1          & -          \\
Linear     & fc2\_4        & Fully Connected & 1          & -          \\
Output     & output     & -               & 4          & -          \\ \hline
\end{tabular}
\caption{Network architecture}
\label{tab:network_architecture}
\end{table}

    Since the dual equation \eqref{discrete_vector} is solved through the information from primal solutions, we trained the model with inputs to be the two-dimensional location information, the Mach number and attack angle's configurations, the volume of a specific element, the four-dimensional solution information and four-dimensional residual information. Then the outputs are the dual solutions.

    In the training process, we adopted the initialization method proposed in \cite{he2015delving}. It provides a better weight initialization for Rectified Linear Unit(ReLU)-based networks, leading to faster convergence. Since the activation functions in our neural networks are ELU, we adopted the leaky ReLU activation function during the initialization process. In order to learn the characteristics of the data better, the inception module is introduced to help CNNs learn more effective and complex patterns. Besides, since the constructed neural networks contain multiple layers, the initialization method addresses the vanishing gradient problem more effectively for networks using ReLU activation functions, ensuring that the gradients are propagated effectively throughout the layers. Moreover, the residual connection method is applied as well. We connect the residual at the end of the encoder and decoder parts respectively so that the training can enable information transmitted in such complex architecture. It also adds regularization to the model, preventing overfitting by providing a form of parameter sharing across layers\cite{he2016deep}. 
    
    Since the real value of dual solutions is mostly lower than 0.01, the general loss function may not be suitable to generate a reliable value. The mean square error is set as loss function at first. Even though the convergence of loss function is smooth, the trained model may not generate reliable values for the four dual variables simultaneously. Then the loss function we considered is the summation of relative errors in the four outputs, i.e.,
    \begin{equation}
      Loss = \sum_i \left|\frac{\tilde{z}_i-z_i}{z_i}\right|.
    \end{equation}
    Here $z_i$ is the target value while $\tilde{z}_i$ is the predicted value. By using Bayesian retrieval methods, we found the Adam optimizer best fits the training process.

\begin{figure}[p]\centering
        \frame{\includegraphics[width=0.48\textwidth]{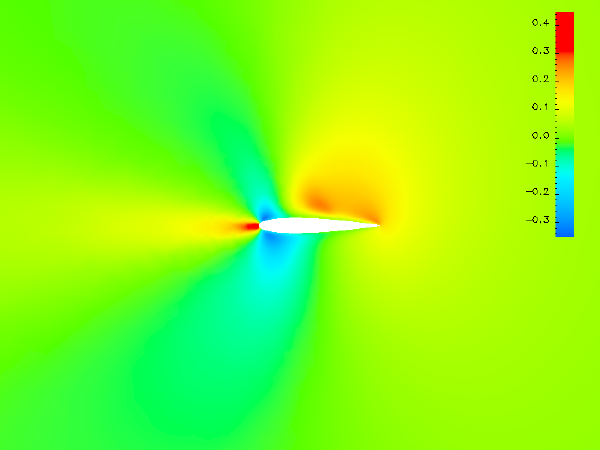}}
        \frame{\includegraphics[width=0.48\textwidth]{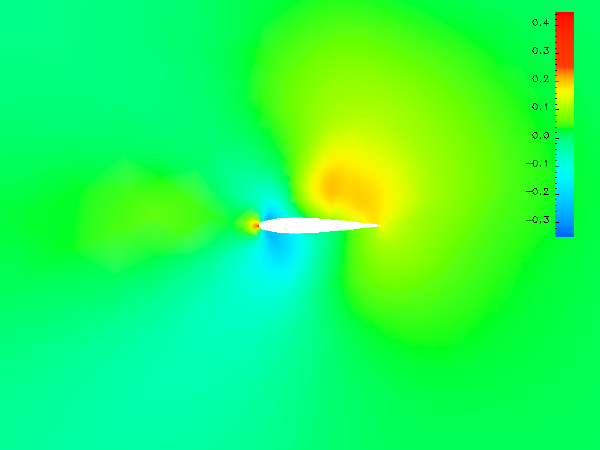}}   \\
        \frame{\includegraphics[width=0.48\textwidth]{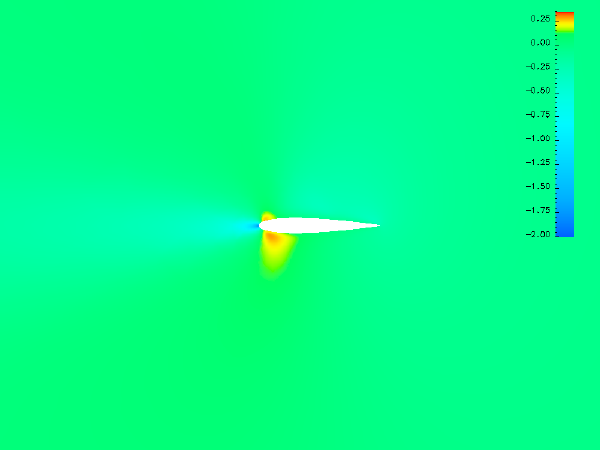}}
        \frame{\includegraphics[width=0.48\textwidth]{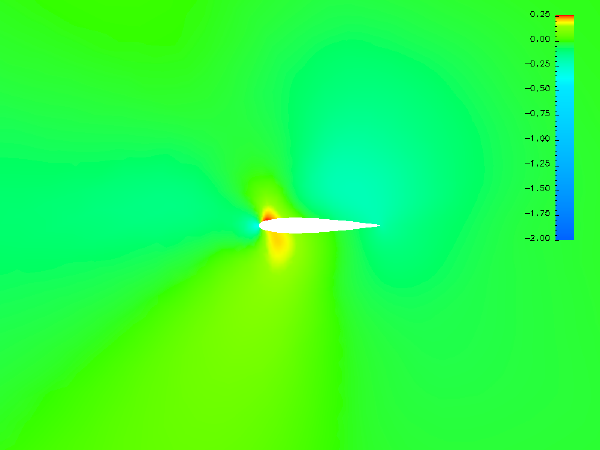}}\\
        \frame{\includegraphics[width=0.48\textwidth]{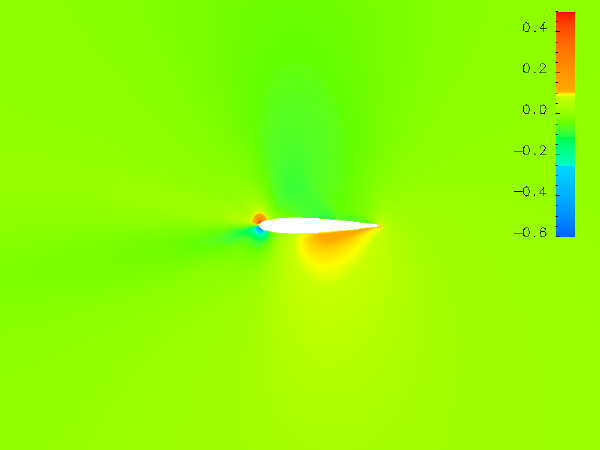}}
        \frame{\includegraphics[width=0.48\textwidth]{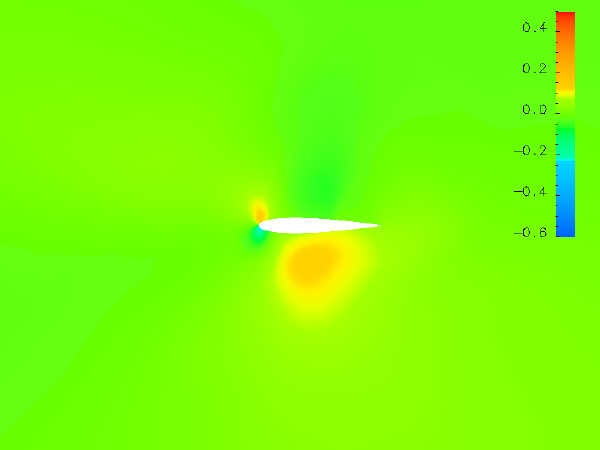}} \\  
        \frame{\includegraphics[width=0.48\textwidth]{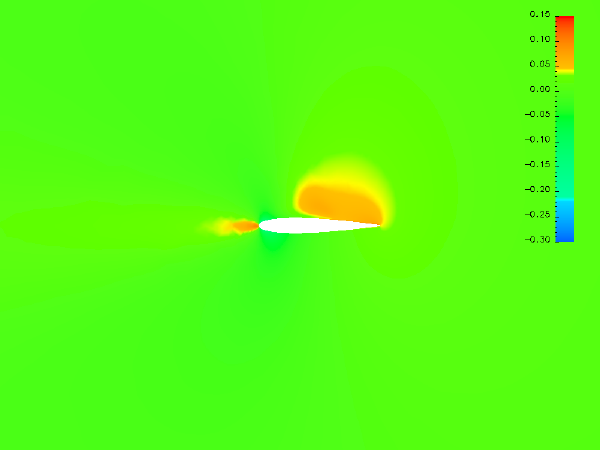}}
        \frame{\includegraphics[width=0.48\textwidth]{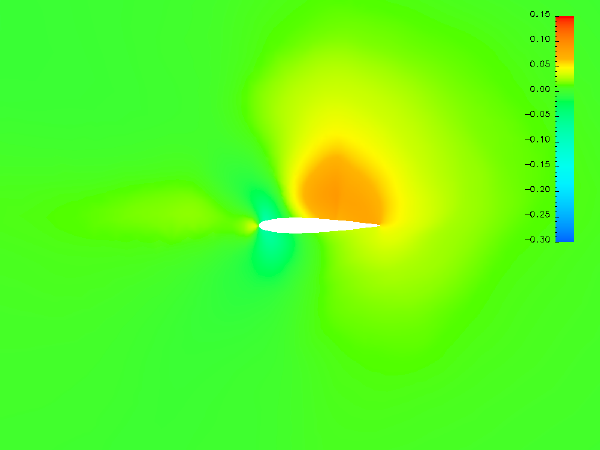}}

      \caption{Left column: The four dual variables generated by GMG solver; Right column: The four dual variables generated by CNNs model; }
      \label{dualCNNsolver}
      \end{figure}

    The data we trained contains the simulation in different configurations. To mitigate the model overfitted on specific models, we adopted the cross-validation technique. It involves partitioning the data into subsets, training the model on some of these subsets, the training set, and validating the model on the remaining validation set. By averaging the model's performance over different subsets, we can reduce the variability and get a better understanding of how the model will perform on unseen data. The result of the training with this cross-validation process can be seen below. The behavior of this model can be seen in Figure \ref{dualCNNsolver}. The figures indicate that the trained model can capture the structures of the dual solutions. Details about the practical calculation will be introduced in the numerical experiment section.
    \begin{figure}[h]\centering
        \includegraphics[width=0.8\textwidth,height=0.4\textheight]{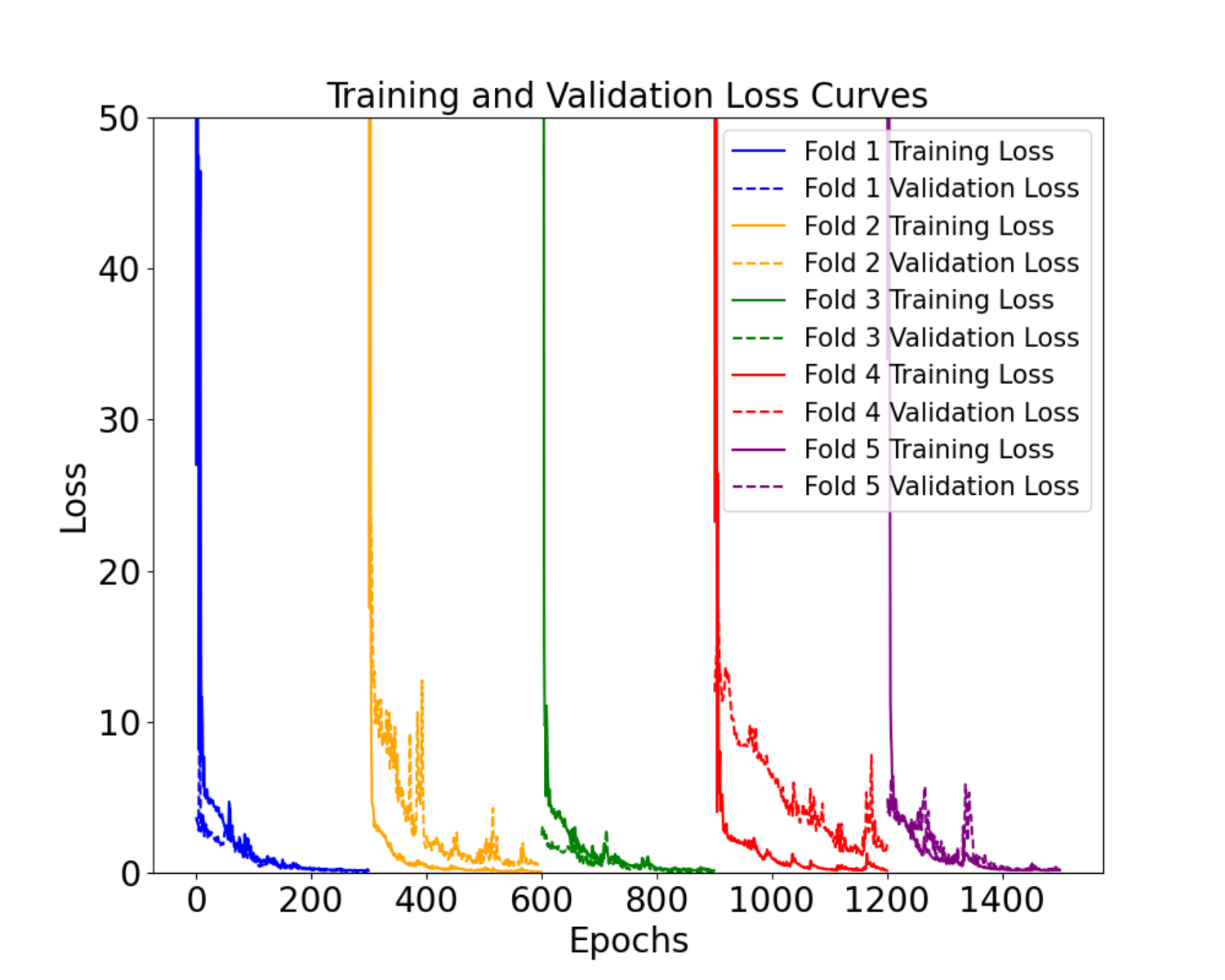}
      \caption{Convergence histograms of the training with cross-validation}  
      \label{resTrace}
      \end{figure}

      The result in Figure \ref{resTrace} is a test on the training for airfoil NACA0012 with $0.8$ Mach number, $1.25^{\circ}$ attack angle. The loss can converge to ideal precision when the iteration is processed. The fold2 may not behave like other examples since the data size is not sufficiently large which will be remitted in a more complex data set.
      \subsection{Dual consistency for training datasets in CNNs}
\begin{figure}[p]\centering
        \frame{\includegraphics[width=0.48\textwidth]{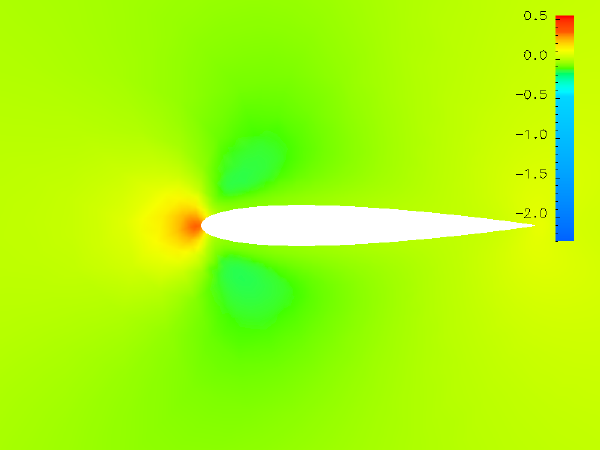}}
        \frame{\includegraphics[width=0.48\textwidth]{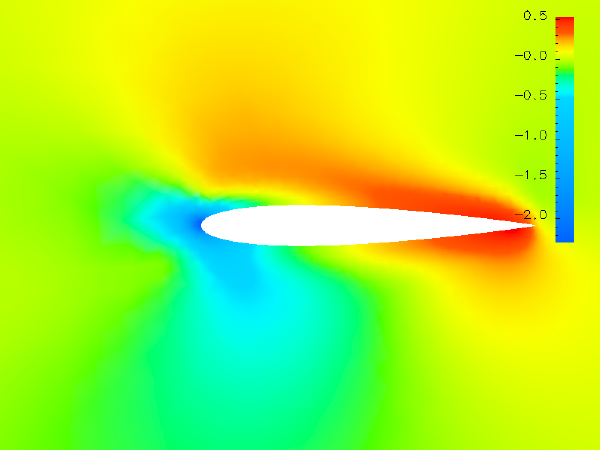}}   \\
        \frame{\includegraphics[width=0.48\textwidth]{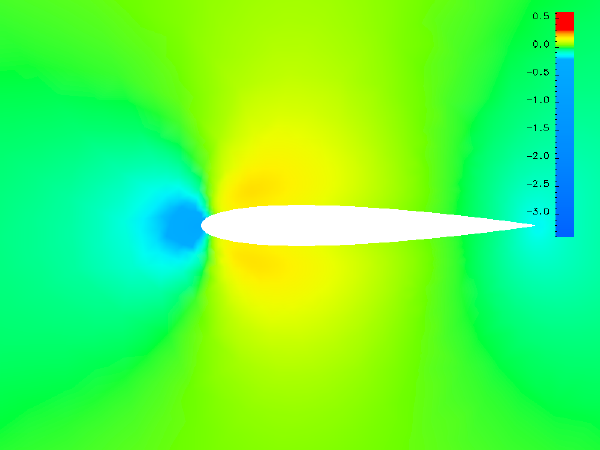}}
        \frame{\includegraphics[width=0.48\textwidth]{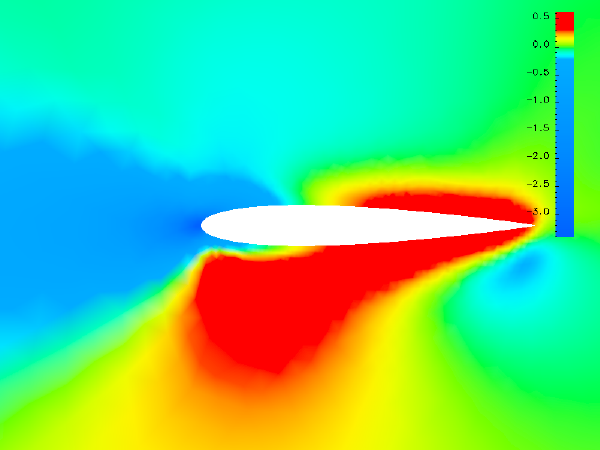}}\\
        \frame{\includegraphics[width=0.48\textwidth]{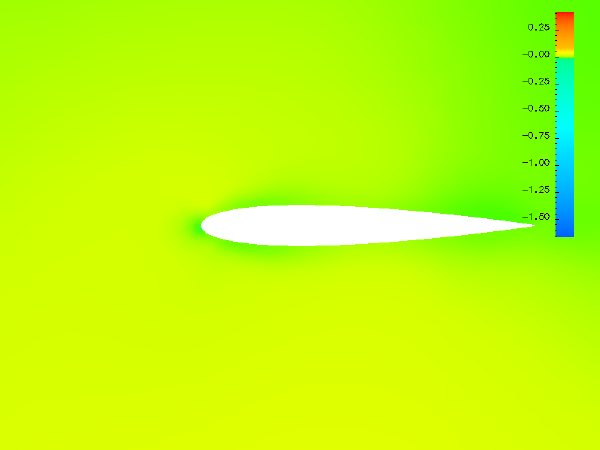}}
        \frame{\includegraphics[width=0.48\textwidth]{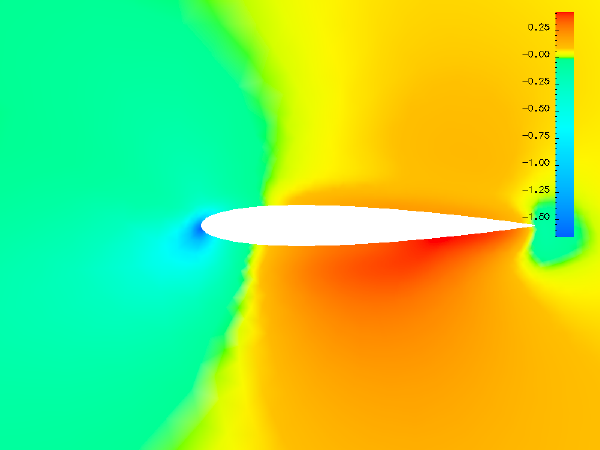}}   \\
        \frame{\includegraphics[width=0.48\textwidth]{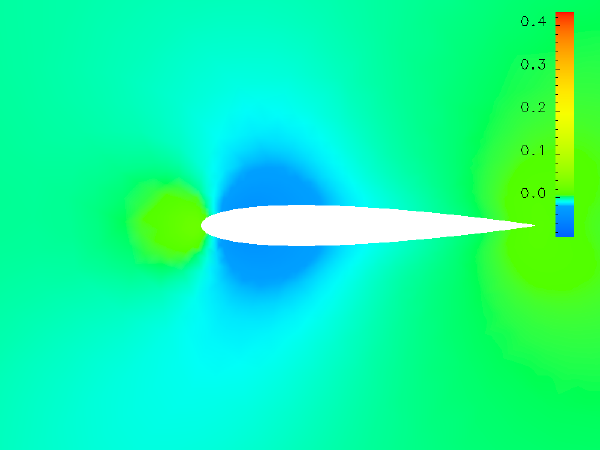}}
        \frame{\includegraphics[width=0.48\textwidth]{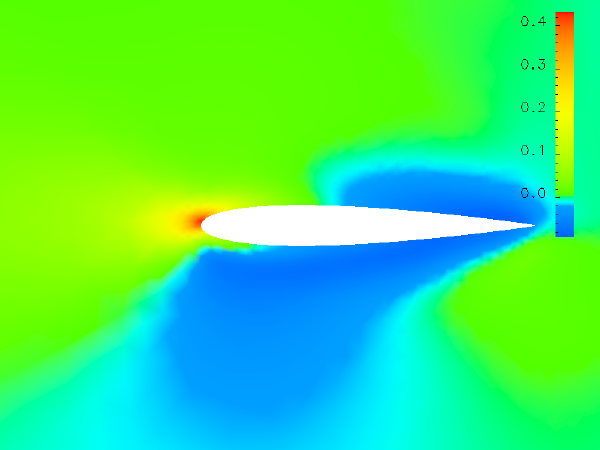}}

      \caption{Left column: The four dual variables generated by CNNs solver trained from dual-consistent datasets; Right column: The four dual variables generated by CNNs solver trained from dual-inconsistent datasets; }
      \label{dualconsistency}
      \end{figure}

      In \cite{wang2023towards}, we analyzed the importance of dual consistency for the DWR-based adaptation in Newton-GMG solver. If the dual consistency is not satisfied, additional waste on the computational resources will occur\cite{hartmann2007adjoint,dolejsi2022}. Worse still, the adaptation may lead to a mesh with target function with lower error. We test the training data from NACA0012 with 0.5 Mach number, 0 attack angle, and zero normal velocity boundary condition. In this study, we validate that if the training datasets are obtained from a dual-inconsistent solver, the performance of CNNs prediction will perform badly as well. The convolutional neural networks can not capture the relation between dual and primal solvers. Then the training process misunderstands the distribution, leading to even worse dual variables. Conversely, in the datasets obtained from a dual-consistent solver, the smooth and symmetric properties can be preserved in this model as shown in Figure \ref{dualconsistency}. Then, it indicated that a dual-consistent DWR-based solver should be constructed at first. Based on that, an efficient CNNs form dual solver can be constructed.

      It is worth noting that the datasets of dual solutions we obtained are from the GMG solver. As the dual-inconsistent scheme may pollute the adaptation around the boundary, some unexpected singularities may occur. Then the dual equations may generate a linear system with a loss of regularization. In order to resolve this issue, a regularization term was added to the dual equations, leading to a GMG solver with a more stable convergence rate. To handle more complicated boundary conditions, the boundary modification techniques proposed by Hartmann\cite{HARTMANN2015754} should be introduced.
      \section{Acceleration Strategy}
      
      The trained model above is saved as the ONNX form so that the C++ library maintained by our group AFVM4CMD can invoke the dual solver. AFVM4CFD is a sophisticated solver that can efficiently handle steady Euler equations. This solver incorporates modules such as k-exact reconstruction, Bessel curves, and geometrical multigrid techniques, thereby offering robust numerical solutions.

      With the introduction of the trained model, which substitutes the conventional dual solver, the DWR-based mesh adaptation process is notably expedited. To further bolster performance and accelerate the simulation process, we implemented additional enhancements within this framework. The subsequent subsections will elaborate on this module that significantly optimizes the numerical experimentation process.
      
\subsection{An Automatic Adjustment of Tolerance in h-adaptive process}
    Even though the DWR-based mesh adaptation method generates high-quality indicators, formulating an appropriate tolerance presents challenges. Different from the general adaptation issues, setting the tolerance too low in DWR-based adaptation can result in subsequent steps resembling a uniform refinement. This disobeys the initial intent of economically solving the quantity of interest. Conversely, if the tolerance is set too high, the adaptation strategy may fail to capture the region which influences the target function significantly. Motivated by the decreasing threshold method proposed in \cite{nemec2007adjoint}, the indicators are plotted for analysis and a straightforward algorithm is proposed to adjust the tolerance dynamically. The indicator we adopted is the multiplication of dual solutions and residuals, i.e.,
    
    \begin{equation}
      \eta_{T_i} = \sum_{T_{ij}}|(\mathbf{z}_h)^T\mathcal{R}_h(\mathbf{u}_h^H)|,
    \end{equation}
    where $T_{ij}$ is the subelements $T_{i}$. We recorded the indicators statistically in three consecutive refined meshes. The distribution concentrated from the left end to the right. Then we tested the distribution function with the Kolmogorov-Smirnov method. The result shows that the distribution may approximate Weibull distribution with KS-stat equals $0.98327$ in this example.
    \begin{figure}[]\centering
    \includegraphics[width=0.32\textwidth]{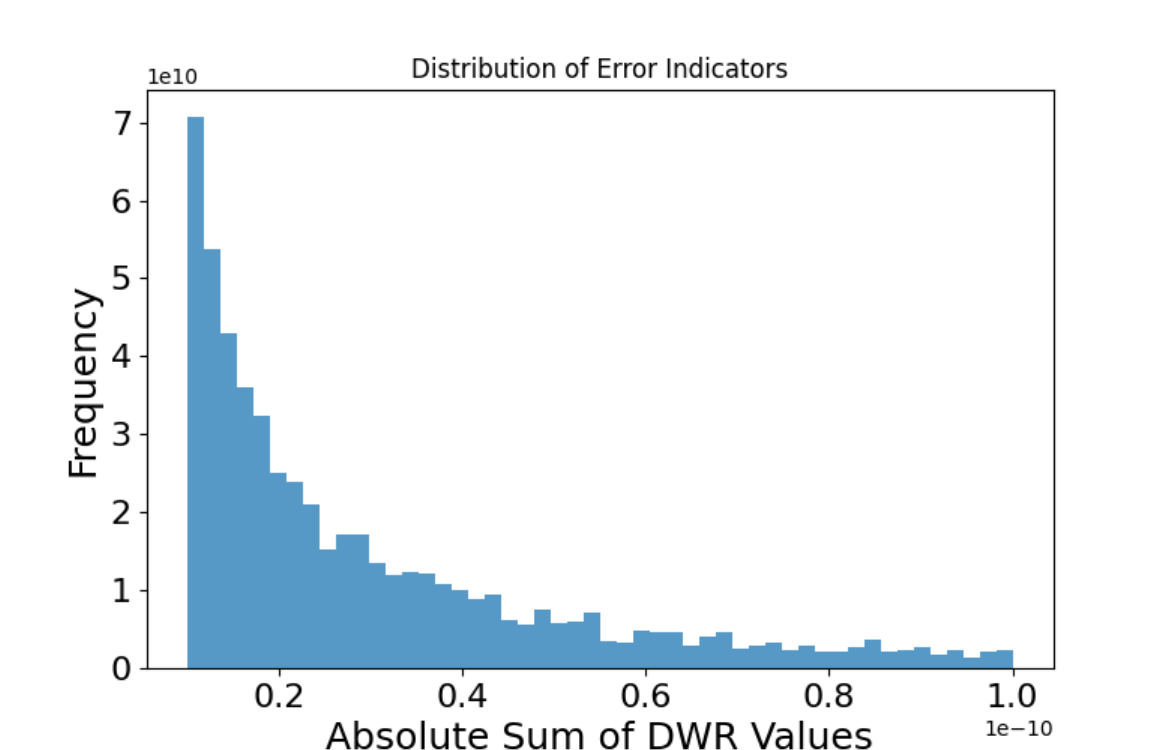}
    \includegraphics[width=0.32\textwidth]{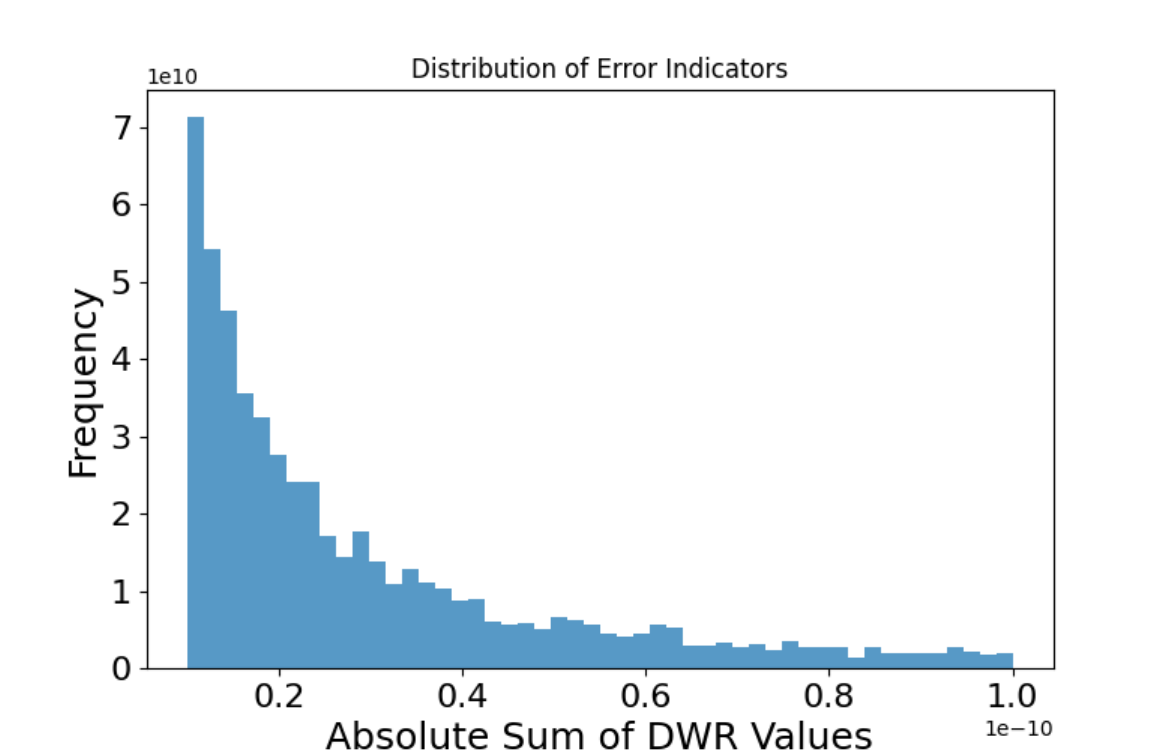}
    \includegraphics[width=0.32\textwidth]{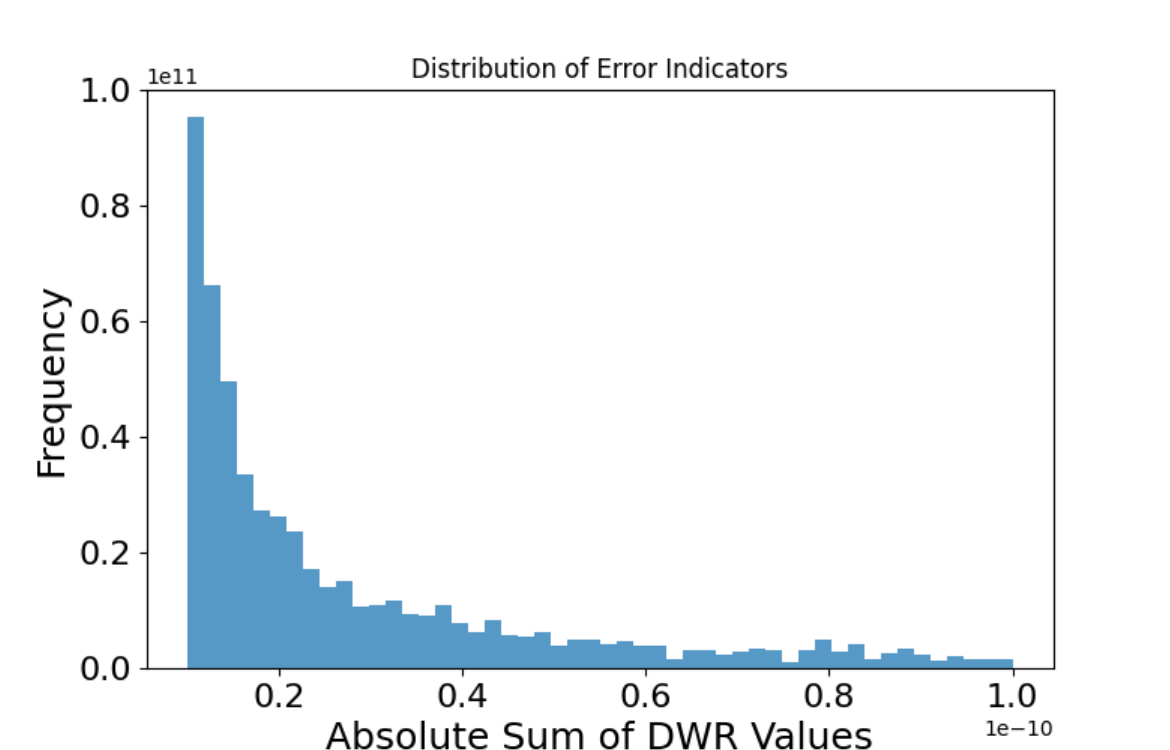}
    \caption{Indicators in three consecutive refined meshes}
    \label{indicator3}
    \end{figure} 

    However, when the configurations changed, the distribution may approximate gamma distribution in some circumstances. Then it's not robust to choose the tolerance with the expression of the distribution function. Alternatively, we designed a well-implemented strategy. The indicators in each element of a given mesh are sorted at first. Then, tolerance is set to be the indicator whose index is a fixed proportion to the mesh size. In this way, with the mesh adaptation processed, the tolerance can be adjusted dynamically.

\begin{figure}[ht]\centering
      \frame{\includegraphics[width=0.31\textwidth]{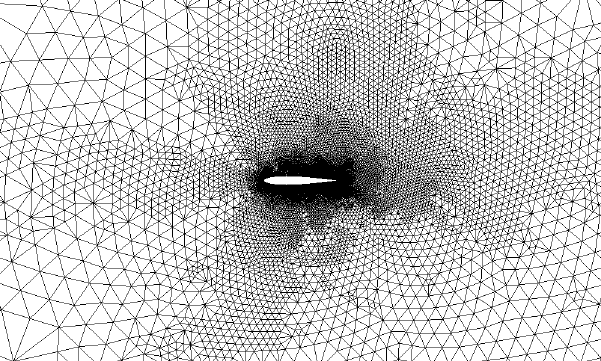}}   
      \frame{\includegraphics[width=0.31\textwidth]{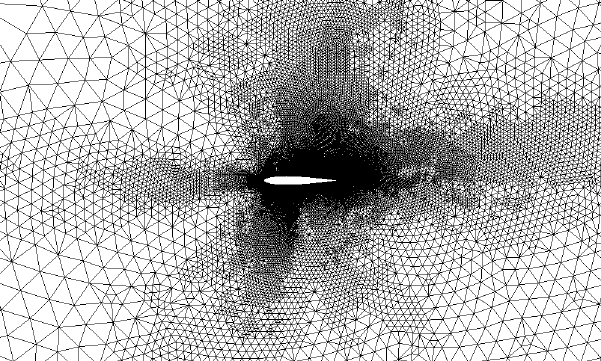}}
      \frame{\includegraphics[width=0.31\textwidth]{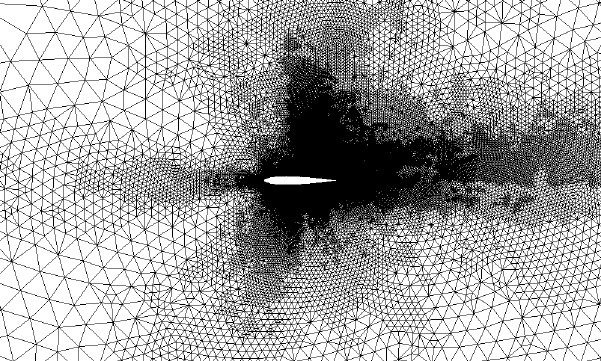}}
    \caption{Consecutive refined meshes with dynamic tolerance strategy. NACA0012, Mach number 0.8, attack angle 1.25$^\circ$.}
    \label{dynamicANDconstant}
    \end{figure}

    To test the dynamic tolerance strategy, numerical experiments are conducted in this part, and we compared the result with four times uniformly refined mesh. After introducing the dynamic tolerance strategy, the number of elements is growing steadily, which means that the corresponding growth of number of elements can be implemented as expected to achieve the required error. From the result in Table \ref{tab:tolerance}, if the constant tolerance is set high, the adaptation may not capture the area which influences the target functional mostly, leading to a rough precision. Conversely, if the constant tolerance is set low, the adaptation at the beginning stage may behave like a uniform refinement, which costs a long time for the calculation of dual equations. However, the dynamic tolerance strategy preserves a stable growth rate, and an expected precision while saving time for generating dual solutions compared with the lower constant tolerance strategy. As shown in Figure \ref{dynamicANDconstant}, the dynamic tolerance strategy performs well for generating the mesh which calculates the target functional precisely.

   \begin{table}[ht]
      \scriptsize
      \centering
      \begin{tabular}{|l|l|c|c|c|c|}
      \hline
      Tolerance& & \multicolumn{4}{c|}{Adaptation} \\
      \cline{3-6}
      & & 1st & 2nd & 3rd & 4th  \\
      \hline
      \multirow{5}{*}{Dynamic} & Number of Elements & 11130 & 33945 & 104985 & 329901 \\
      \cline{2-6}
      & Refined Rate & 0.69 & 0.68 & 0.70 & 0.71\\
      \cline{2-6}
      & Tolerance Value & $4.70936\times 10^{-11}$ & $2.24686\times 10^{-11}$ & $9.77698\times 10^{-12}$ & $1.50204\times 10^{-12}$ \\
      \cline{2-6}
      & Time for Dual& 28 & 111 & 420 & 1811  \\
      \cline{2-6}
      & Error & 0.0132936& 0.0030573 & 0.0007455 & 0.0001009 \\
      \hline
      \noalign{\vspace{3mm}}
      \hline
      \multirow{4}{*}{\begin{tabular}{@{}c@{}}Constant 1\\ $3.0\times 10^{-11}$\end{tabular}} & Number of Elements & 11961 & 31581 & 75363 & 148020 \\
      \cline{2-6}
      & Refined Rate & 0.76 & 0.55 & 0.46 & 0.32  \\
      \cline{2-6}
      & Wall Time for Dual (seconds)& 28 & 120 & 377 & 1186  \\
      \cline{2-6}
      & Error & 0.0133177 & 0.003120 & 0.0015794 & 0.0013368 \\
      \hline
       \noalign{\vspace{3mm}}
       \hline
      \multirow{4}{*}{\begin{tabular}{@{}c@{}}Constant 1\\ $3.0\times 10^{-12}$\end{tabular}} & Number of Elements & 14529 & 50781 & 134772 & 298170  \\
      \cline{2-6}
      & Refined Rate & 0.99 & 0.83 & 0.55 & 0.40  \\
      \cline{2-6}
      & Wall Time for Dual (seconds)& 28 & 150 & 664 & 2442  \\
      \cline{2-6}
      & Error & 0.0133227 & 0.0029431 & 0.0006927 & 0.0002414 \\
      \hline
      \end{tabular}
      \caption{Comparison of Dynamic and Constant Tolerance Strategies}
      \label{tab:tolerance}
    \end{table}

      From the above post-processing technique, we developed a more efficient framework for solving the quantity of interest with the CNNs solver. The algorithm is organized as follows:

      \begin{algorithm}[H]
        \SetAlgoLined
        \KwData{Initial $\mathcal{K}_H$, $TOL$}
        \KwResult{$\mathcal{K}_h$}
        Using the Newton-GMG to solve $\mathcal{R}_{H}(\mathbf{u}_H)=0$ with residual tolerance $1.0\times 10^{-3}$\;
        Convert the element information to Tensors' form\;
        Generating the dual solutions with the trained ONNX model\;
        Calculate the error indicator for each element\;
        Select the dynamic tolerance for adaptation\;
        \While{$\mathcal{E}_{K_H}>TOL$ for some $K_H$}{
          Adaptively refine the mesh $\mathcal{K}_{H}$ with the process in \cite{HU2016235};}
        \caption{DWR for one-step mesh refinement}
      \end{algorithm}

    \subsection{OpenMP parallelization of the algorithm}
Parallel computing enables the simultaneous execution of computations, leading to significant reductions in time, especially for the CFD domain where simulations can be extremely time-consuming. Among these, OpenMP is a widely used parallel programming model that greatly facilitates the execution of tasks by dividing them into multiple threads that can run concurrently. However, it requires a careful design of the algorithm to avoid potential pitfalls. In this part, we mainly concern the scalability for different parts of the algorithm. 

The primal solver, Newton-GMG solver, contains reconstruction, update cell average, and geometrical multigrid module for each Newton iteration step. If the reconstruction patch should be built in each Newton iteration step, it will consume too much time. Then we build a cache module to reserve the information of the reconstruction patch in advance so that the identical information will not be calculated redundantly. Since the update of the cell average is independent for different elements, it exhibits a satisfactory performance of parallelization. 

\begin{figure}[!h]\centering
  \includegraphics[width=0.45\textwidth]{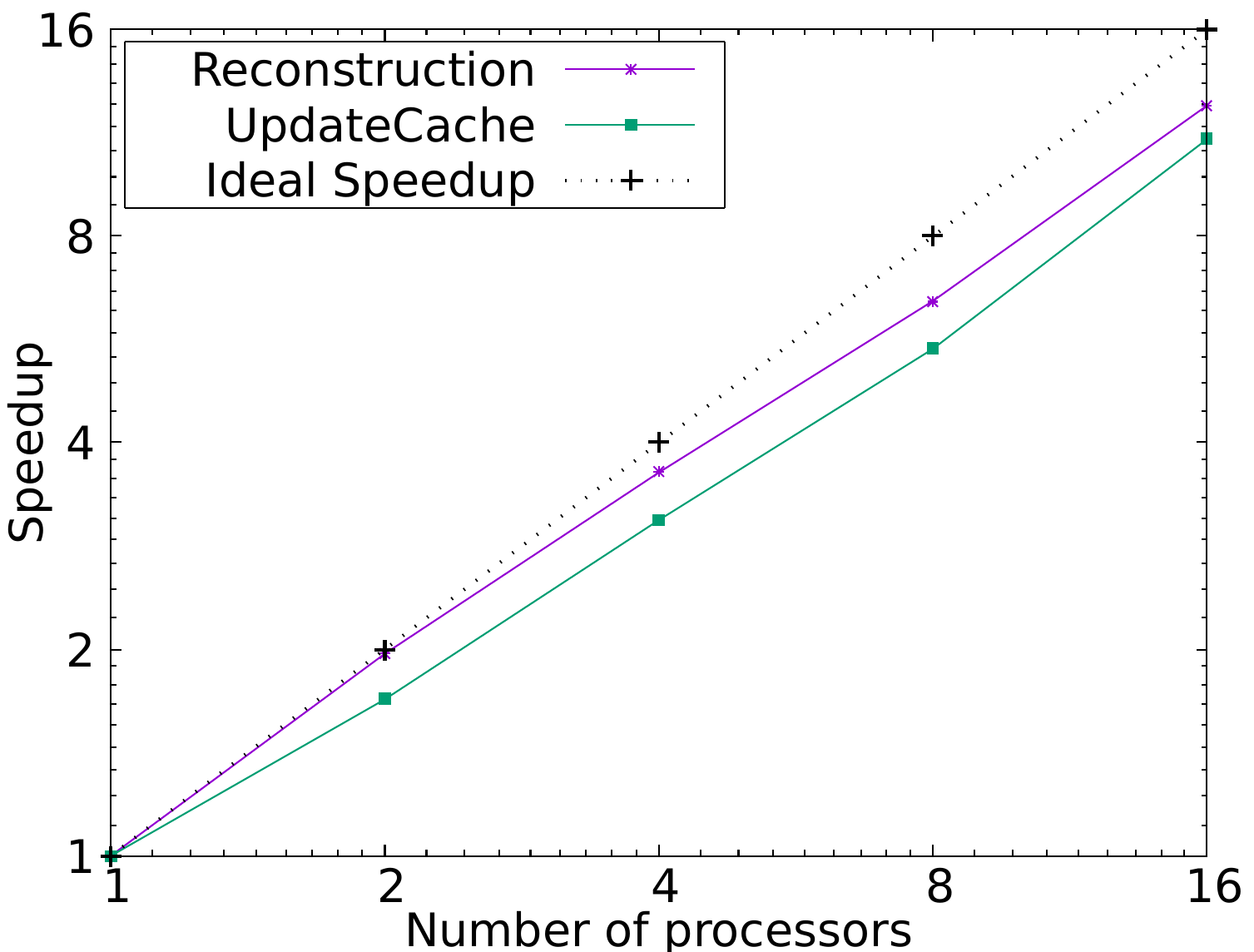}
  \includegraphics[width=0.45\textwidth]{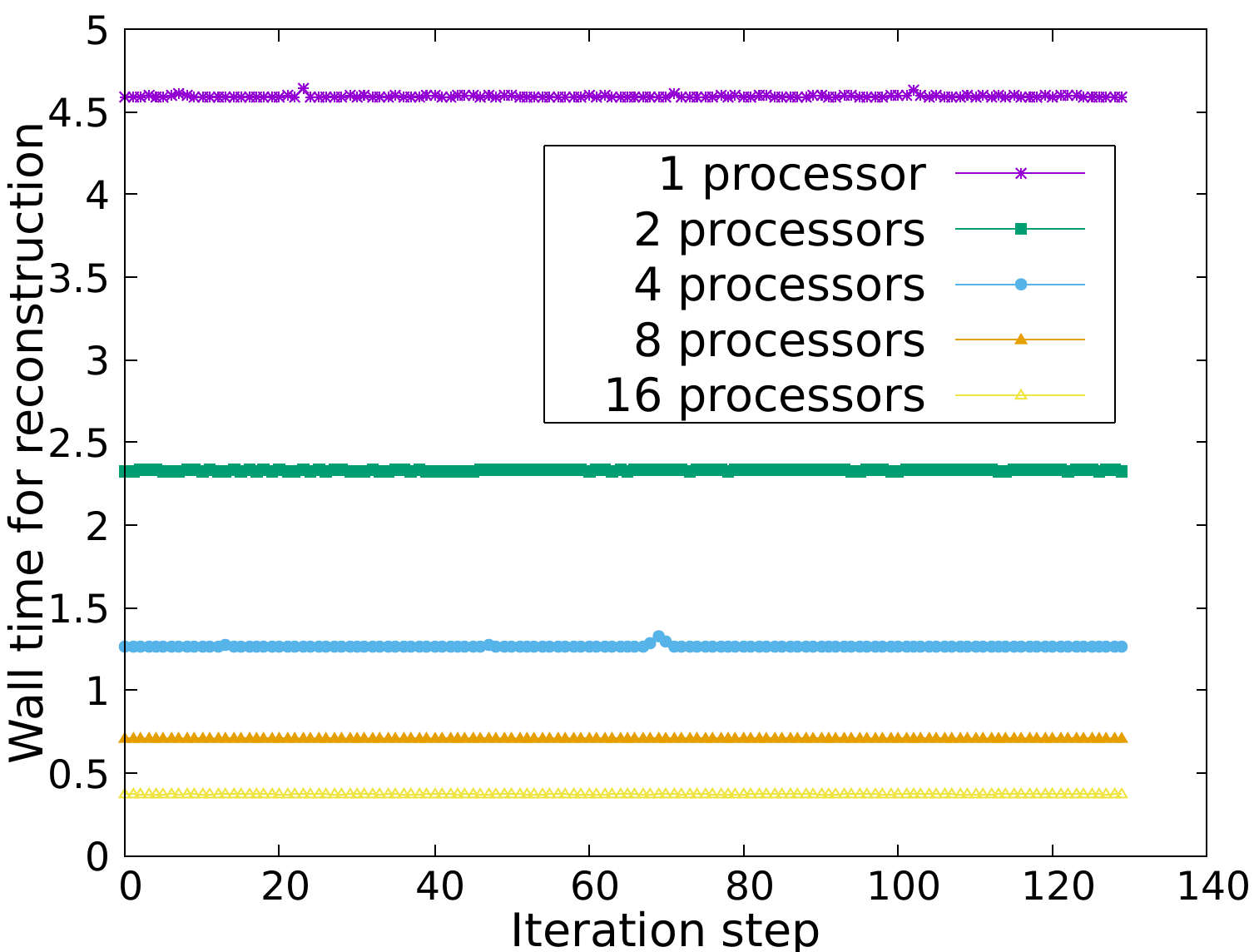}
\caption{Left: Mean time spent on modules during the Newton iteration step with the different number of threads. Right: Time spent on the reconstruction during the whole Newton iteration process. }
\label{Threadcompare}
\end{figure} 
\begin{figure}[!h]\centering
  \includegraphics[width=0.45\textwidth]{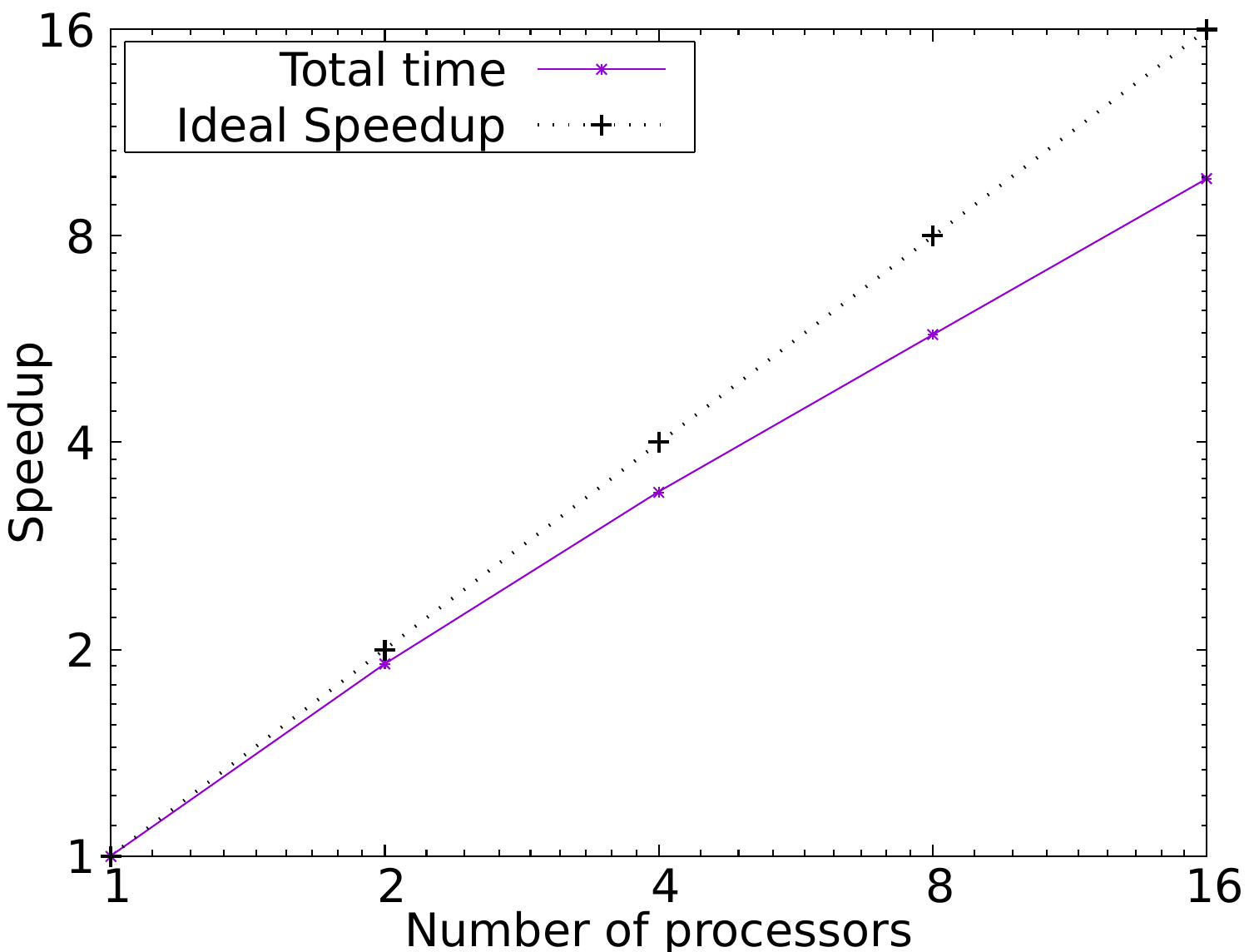}
  \includegraphics[width=0.45\textwidth]{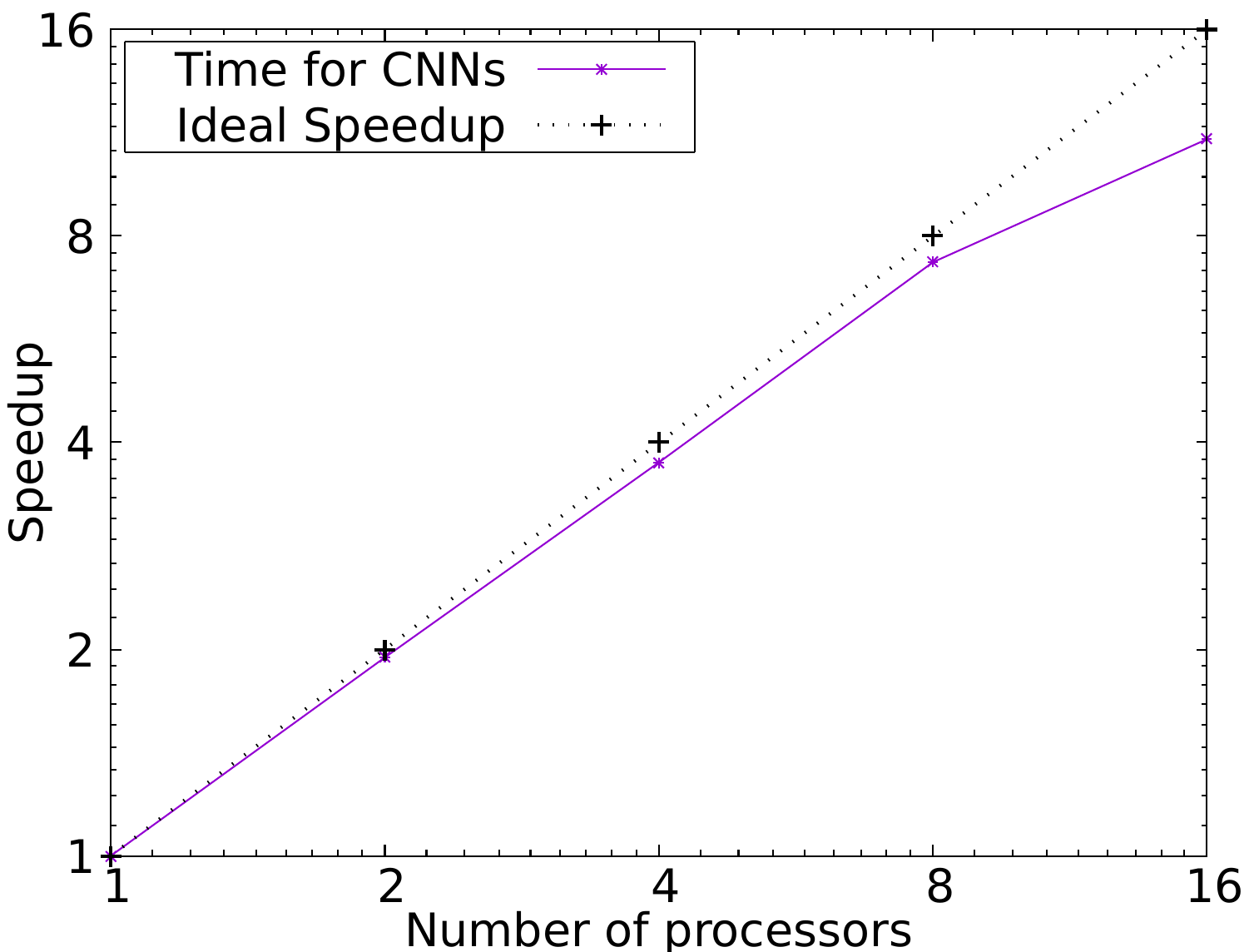}
\caption{Left: Total time spent for the primal solver with the different number of threads. Right: Total time spent on the dual solver with the different number of threads. }
\label{totalThreadcompare}
\end{figure}

The speedup in parallel computing is generally defined as
\begin{equation}
    Speedup = \displaystyle\frac{t(1)}{t(N)},
\end{equation}
where $t(1)$ is the time for running the algorithm with one processor and $t(N)$ for $N$ processors respectively.

Theoretically, the time used for parallel computing will be halved as the number of threads doubles. Then the ideal value of speedup should be $N$ when there exist $N$ processors. However, the performance of parallel computing may not behave as expected due to the complexity of the algorithm. To compare the efficiency of the parallelization in our framework, the problem size remains constant in this section. The experiments in this part are all conducted on RAE2822 airfoil whose Mach number is $0.729$, attack angle $2.31^\circ$ on the mesh with $232,704$ elements.

The performance can be seen in Figure \ref{Threadcompare}. It illustrates that the reconstruction part and update part are well suited for parallelization. With the increase in the number of threads, time spent on the calculation has been saved significantly. However, the lower-upper symmetric Gauss-Seidel iteration in each Newton iteration step is not suitable for parallelization, which cannot be optimized even if more threads are used. Then the scalability for the whole process of the primal solver is limited. Figure \ref{totalThreadcompare} illustrates the effect of scalability of primal solver and dual solver. It shows that parallel computing has significantly enhanced computational efficiency, and the scalability of the dual solver further demonstrates that the method of using neural networks to solve dual solutions can greatly accelerate the calculation. The initialization of the trained ONNX model influenced the scalability and will be improved upon in future work. Even though, parallelization is still an effective technique for accelerating the simulation of our algorithm. Then the following numerical experiments are all conducted with the parallism.

\section{Numerical Experiments}
In this section, we will show the performance of CNNs for generating dual solutions with significant time savings. Then the generalizations of the model are demonstrated in three different aspects. Firstly, even when we deal with the adaptation with data not included in the training sets, the behavior still preserves well till the expected precision. Secondly, the model is trained with different Mach numbers and attack angles. Then it can simulate generalized configurations. Thirdly, since the framework is constructed on an unstructured mesh, it is very important to show that the model works well for a more complicated geometry.

\subsection{Time consumption of CNNs dual solver}
The most important reason we adopt the CNNs solver is its significant acceleration in generating dual solutions. The trained models are tested on the configurations shown in Table \ref{timeSaving}. In each adaptation step, the GMG solver, and CNNs dual solver generate high-quality meshes which can derive the quantity of interest with comparable precision.
The data reveals that the CNNs dual solver can dramatically reduce computational time, particularly as the adaptive process progressively advances. In the Three-Airfoils model, which contains much more complicated geometries, the CNNs solver shows great potential for time saving.

\begin{table}[h]\scriptsize
\centering
\begin{tabular}{|l|c|c|c|c|c|c|c|c|}
\hline
\multirow{3}{*}{Configurations} & \multicolumn{8}{c|}{Wall time for dual in different adaptation step (seconds)}\\
\cline{2-9}
 & \multicolumn{2}{c|}{1st} & \multicolumn{2}{c|}{2nd} & \multicolumn{2}{c|}{3rd} & \multicolumn{2}{c|}{4th} \\
\cline{2-9}
 & CNNs & GMG & CNNs & GMG & CNNs & GMG & CNNs & GMG \\
\hline
NACA0012, 0.8 Mach, 1.25$^\circ$ attack angle  & 9 & 28 & 25 & 111  & 81 & 420 & 295 & 1811 \\ 
\hline
NACA0012, 0.76 Mach, 1.05$^\circ$attack angle & 9 & 29 & 28 & 114 & 93 & 373 & 382 & 2253 \\ 
\hline
NACA0012, 0.82 Mach, -1.65$^\circ$attack angle & 9 & 31 & 27 & 125 & 95 & 519 & 387 & 2708 \\ 
\hline
RAE2822, 0.729 Mach, 2.31$^\circ$attack angle & 9 & 32 & 25 & 113 & 78 & 411 & 258 & 1544 \\ 
\hline
3-NACA0012, 0.8 Mach, 1.25$^\circ$ attack angle & 422 & 2324 & 1145 & 9286 & 6521 & 64613 & - & - \\ 
\hline
\end{tabular}
\caption{Comparision of CNNs solver and GMG solver with wall time for dual in different adaptation steps.}
\label{timeSaving}
\end{table}
\subsection{Generalization and model performance}
\subsubsection{Generalization on the training set}
\begin{figure}[ht]\centering
  \frame{\includegraphics[width=0.31\textwidth]{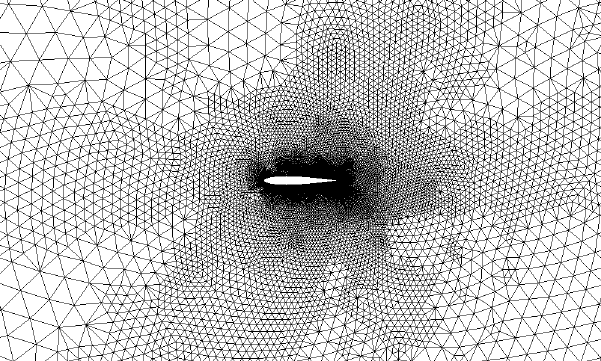}}   
  \frame{\includegraphics[width=0.31\textwidth]{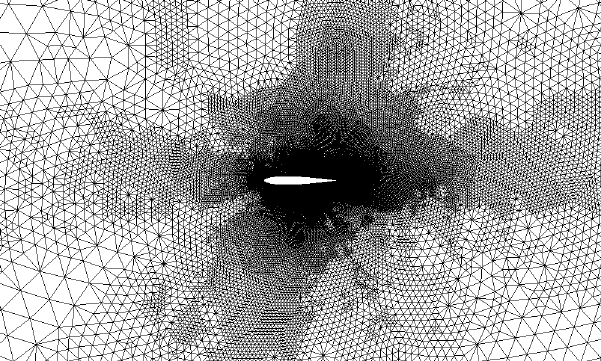}}
  \frame{\includegraphics[width=0.31\textwidth]{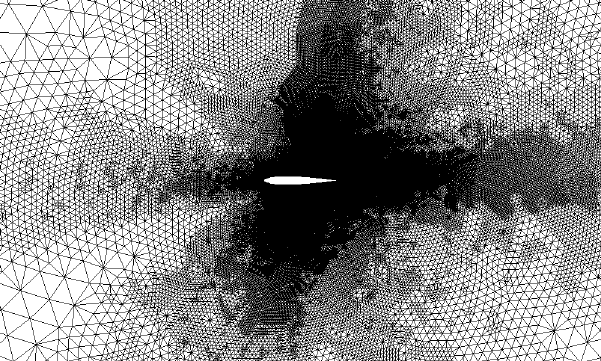}}
\caption{Consecutive refined meshes with CNNs dual solver. NACA0012, Mach number 0.8, attack angle 1.25$^\circ$.  }
\label{CNNadapt}
\end{figure} 

\begin{figure}[!h]\centering
  \includegraphics[width=0.45\textwidth]{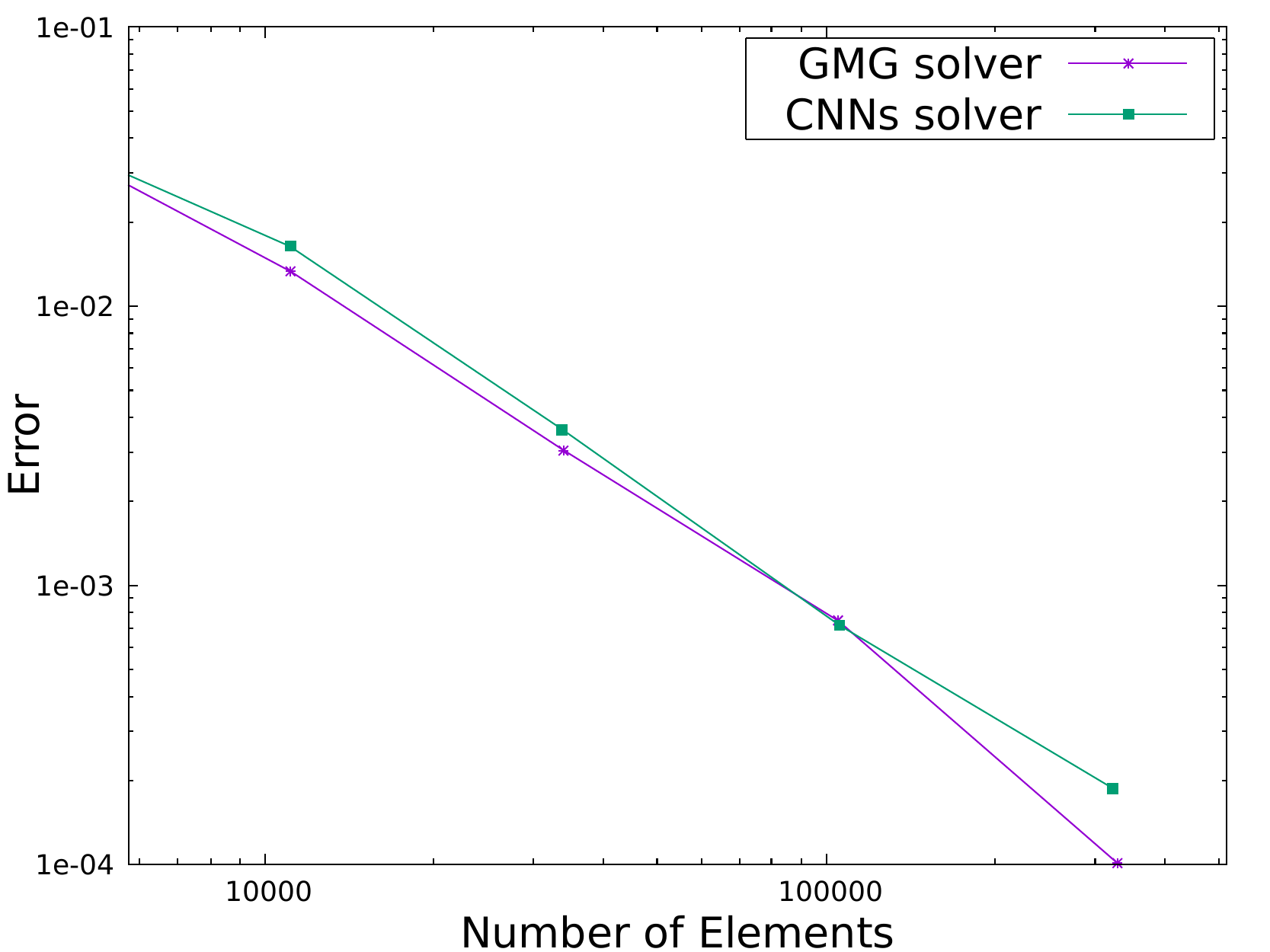}
  \includegraphics[width=0.45\textwidth]{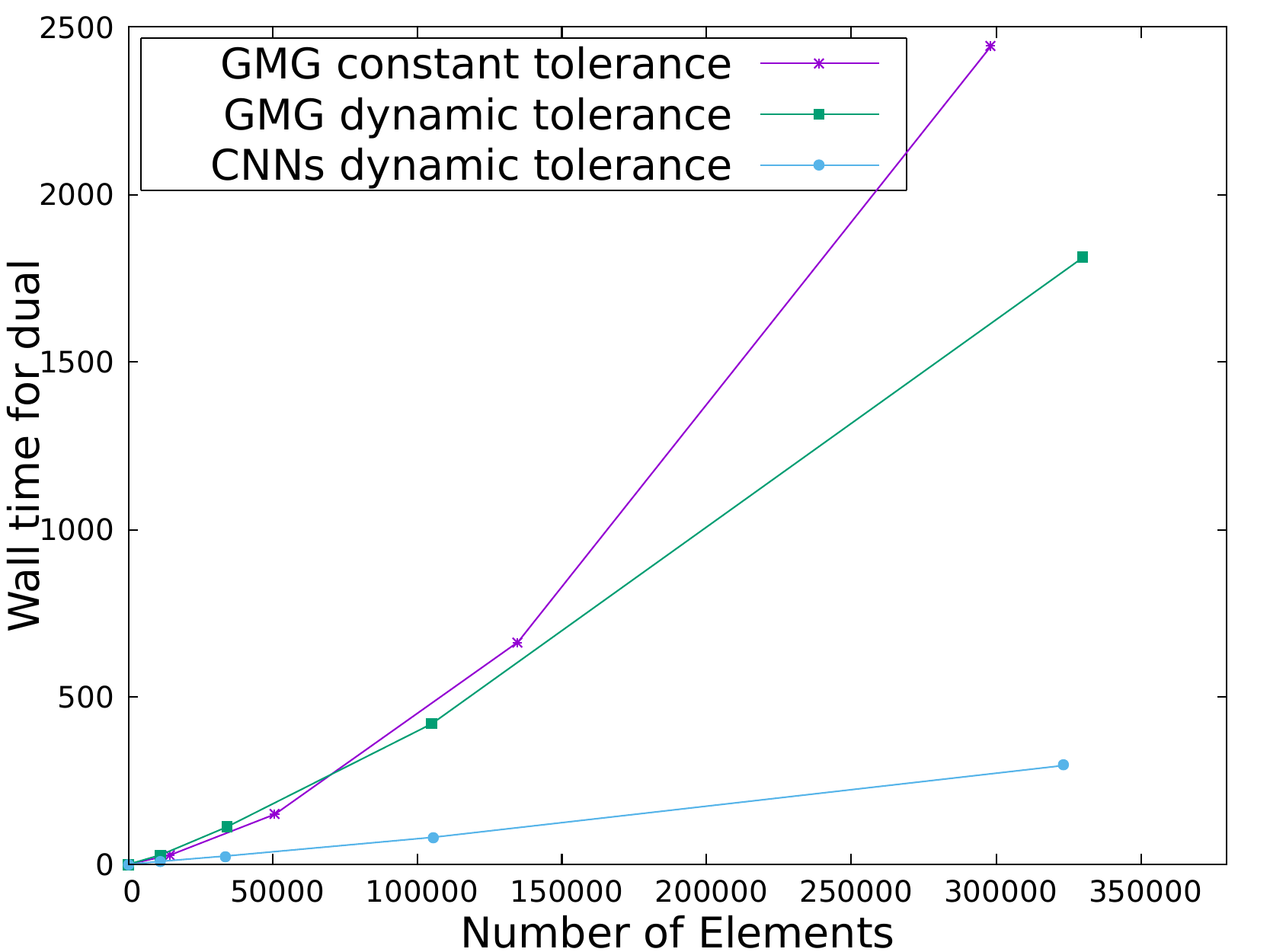}
\caption{Comparison of error and time cost with the dual solver of GMG form and CNNs form. NACA0012, Mach number 0.8, attack angle 1.25$^\circ$. }
\label{DUALvsGMG}
\end{figure} 

In Figure \ref{dualCNNsolver}, the CNNs captured the main structure of dual solutions. With the substituted CNNs form the dual solver, the adaptation gets processed steadily. In Figure \ref{CNNadapt}, it is shown that the dual solver of CNNs form can generate a mesh that balances the dual equations and residuals similar to the GMG solver. The result can be seen from Figure \ref{DUALvsGMG} that the error of CNNs dual solver is comparable to the GMG solver. Here we use the mesh with 4 times uniform refinement with $930,816$ elements to generate a reference value, which is $0.026201$. However, comparing the time for generating the refined mesh, the CNNs form can save time by an order of magnitude. The results in Figure \ref{DUALvsGMG} indicate that time spent on the GMG solver grows exponentially, while the CNNs solver grows linearly. This can be explained that the CNNs solver calls the trained model element by element, then the time complexity of the CNNs solver is just $O(n)$, where $n$ is the number of elements. It is worth noting that to train the CNNs dual solver, we only use the datasets from the mesh with one and two times uniform refinement. In other words, the third and fourth-time refinement is the prediction beyond the training set. As a result, the trained CNNs solver is capable of generating credible dual solutions in this framework.

\subsubsection{Generalization on the configurations}
      \begin{figure}[!p]\centering
        \frame{\includegraphics[width=0.48\textwidth]{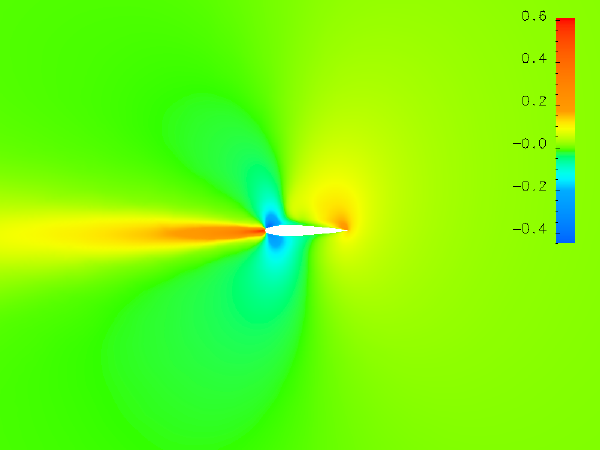}}
        \frame{\includegraphics[width=0.48\textwidth]{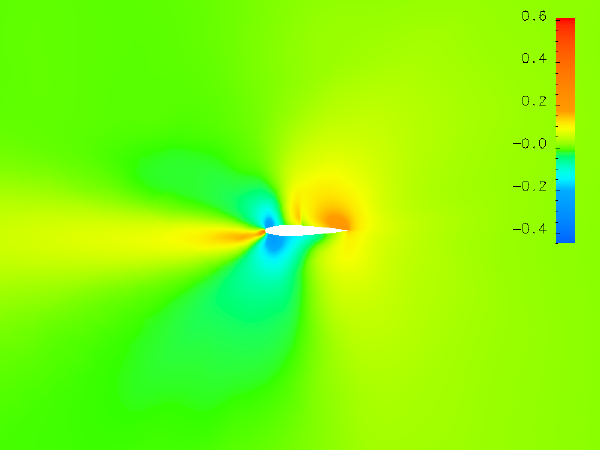}}\\   
        \frame{\includegraphics[width=0.48\textwidth]{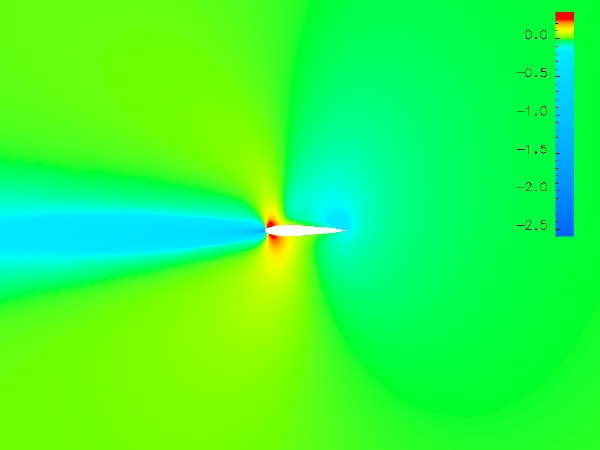}}
        \frame{\includegraphics[width=0.48\textwidth]{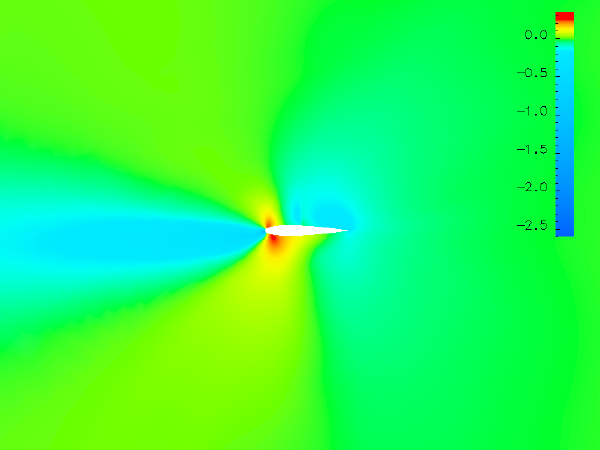}}\\
        \frame{\includegraphics[width=0.48\textwidth]{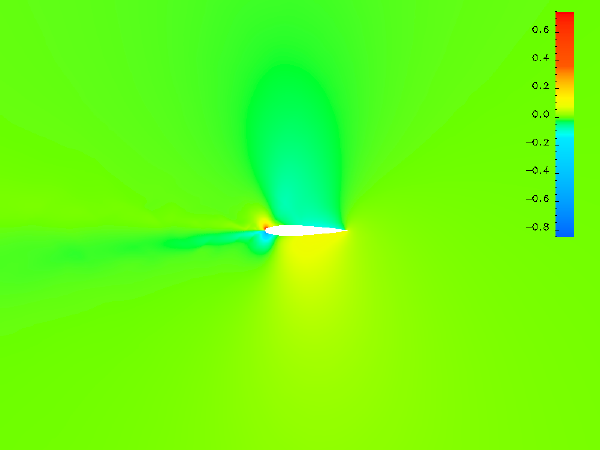}}
        \frame{\includegraphics[width=0.48\textwidth]{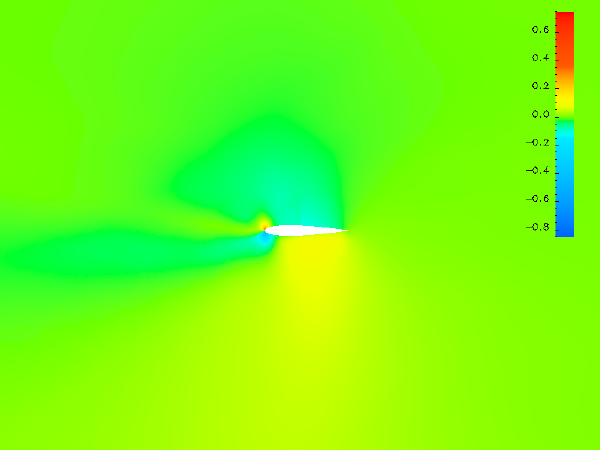}}\\   
        \frame{\includegraphics[width=0.48\textwidth]{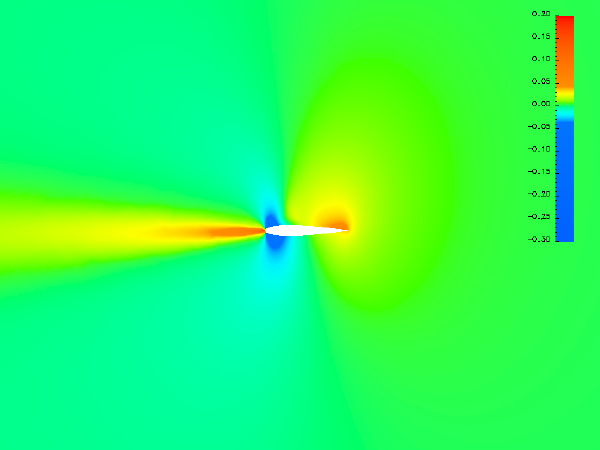}}
        \frame{\includegraphics[width=0.48\textwidth]{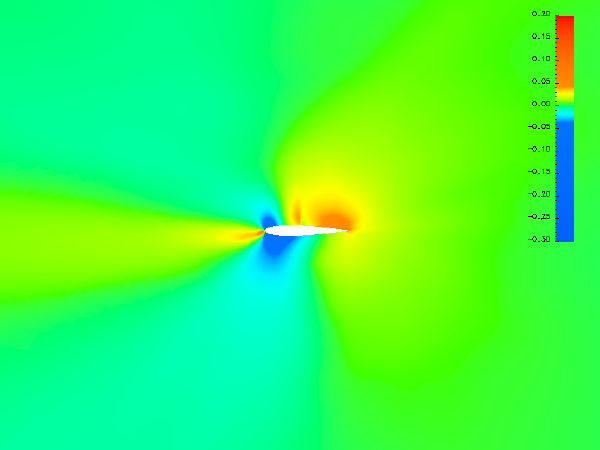}}

      \caption{Left column: The four dual variables generated by GMG solver; Right column: The four dual variables generated by CNNs model; NACA0012, Mach number 0.76, attack angle 1.05.}
      \label{dualCNN0.76solver}
      \end{figure}  
      \begin{figure}[!h]\centering
  \includegraphics[width=0.45\textwidth]{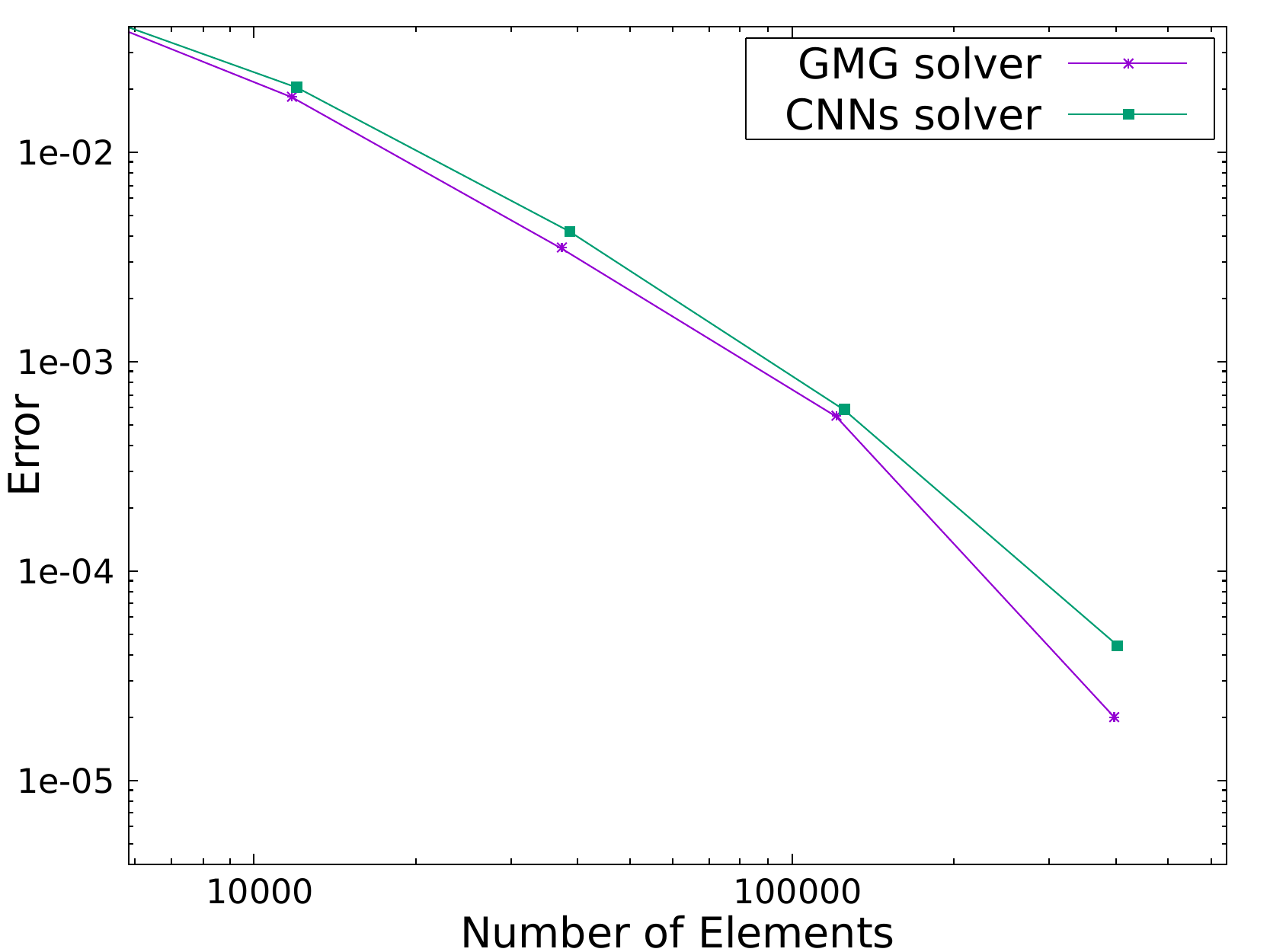}
  \includegraphics[width=0.45\textwidth]{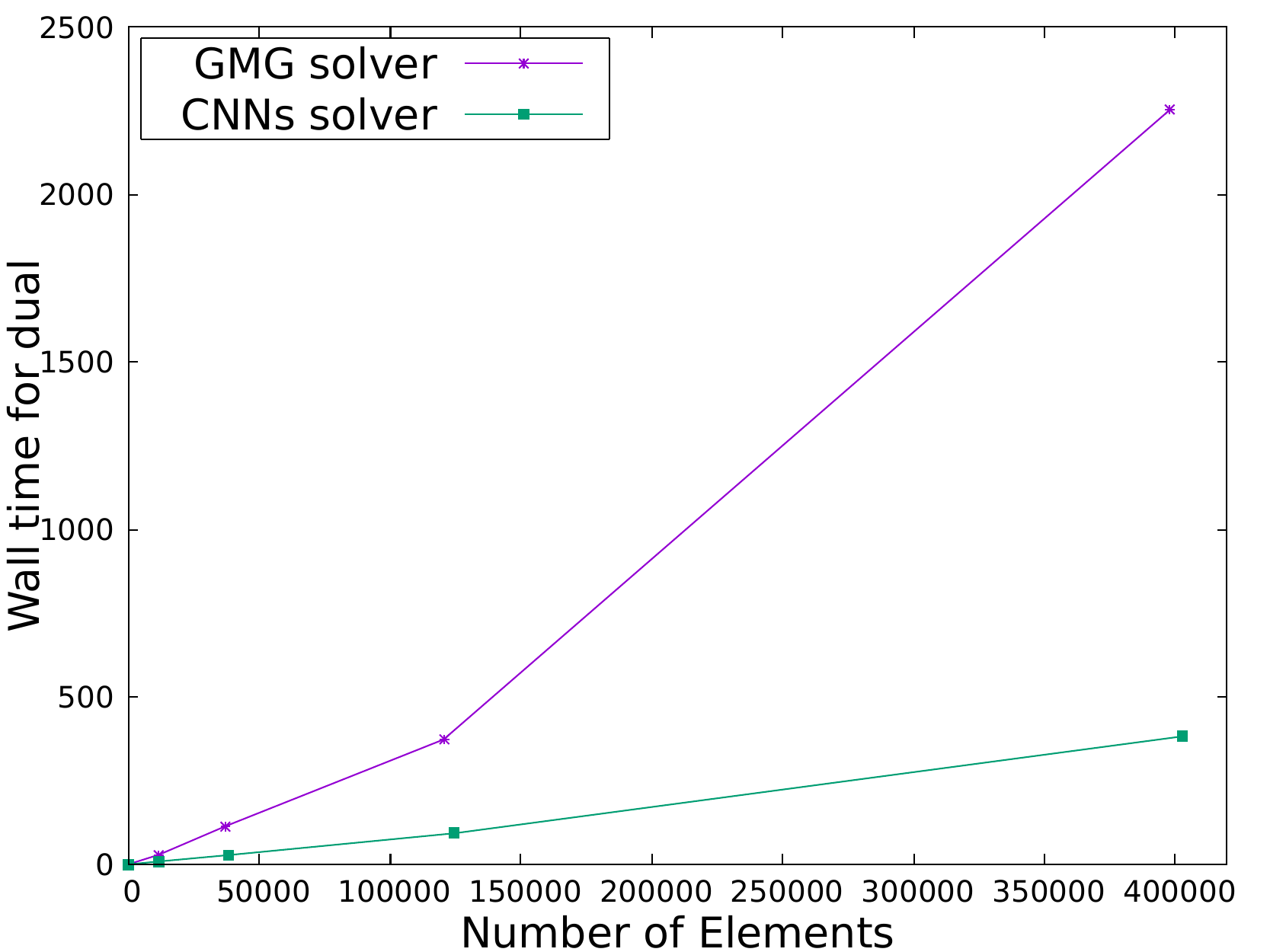}
\caption{Comparison of error and time cost with the dual solver of GMG form and CNNs form. NACA0012, Mach number 0.76, attack angle 1.05. }
\label{DUALvsGMG076}
\end{figure} 
In a shape-optimal design problem, the configurations should be changed for numerous experiments. If the CNNs model needs to be trained every time, the time spent on the training is still expensive. Nonetheless, since the configurations can be treated as input to the CNNs, we can train the model with various configurations so that it can handle different simulations. The CNNs architecture we built is capable of extracting the features between different Mach numbers and attack angles. 
We tested the trained model with the following configurations:
\begin{itemize}
    \item NACA0012, 0.76 Mach, 1.05$^\circ$ attack angle ;
    \item NACA0012, 0.82 Mach, -1.65$^\circ$ attack angle ;
    \item RAE2822, 0.729 Mach, 2.31$^\circ$ attack angle. 
\end{itemize}

    \begin{figure}[!h]\centering
        \includegraphics[width=0.8\textwidth,height=0.4\textheight]{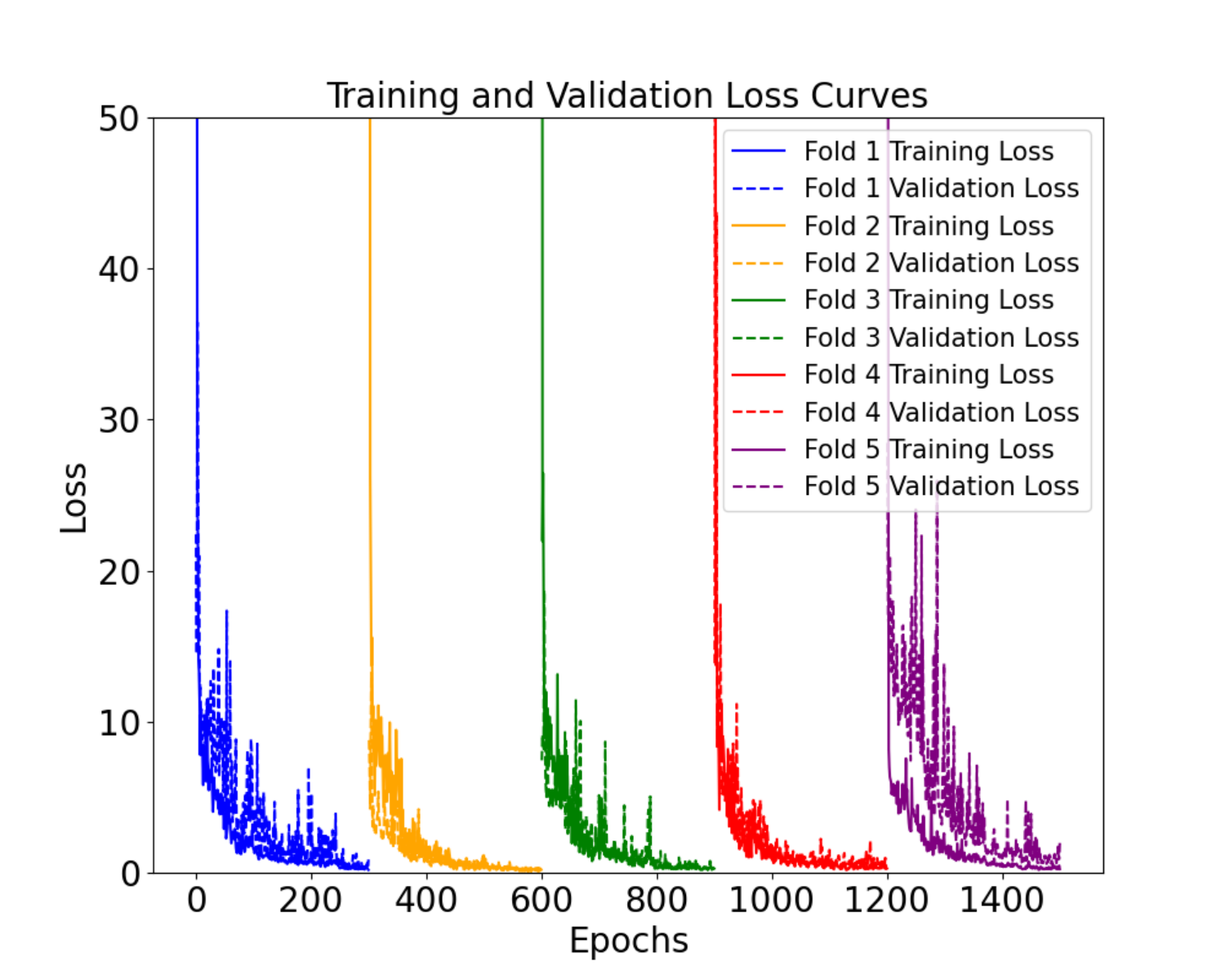}
      \caption{Convergence histograms of the training with different configurations}  
      \label{MultiResTrace}
      \end{figure}
Figure \ref{MultiResTrace} indicates the training process still works well with the cross-validation method.
Then we test the experiment whose configurations are not included in the training sets. Firstly, the Mach number is set at 0.76 while the attack angle is set at 1.05. The distribution of dual variables can be seen in Figure \ref{dualCNN0.76solver}. Similarly, as seen in Figure \ref{DUALvsGMG076}, the CNNs solver generates mesh which solved the target functional in a comparable error with the GMG solver. However, the CNNs solver saves time sharply in this experiment.

      \begin{figure}[p]\centering
        \frame{\includegraphics[width=0.48\textwidth]{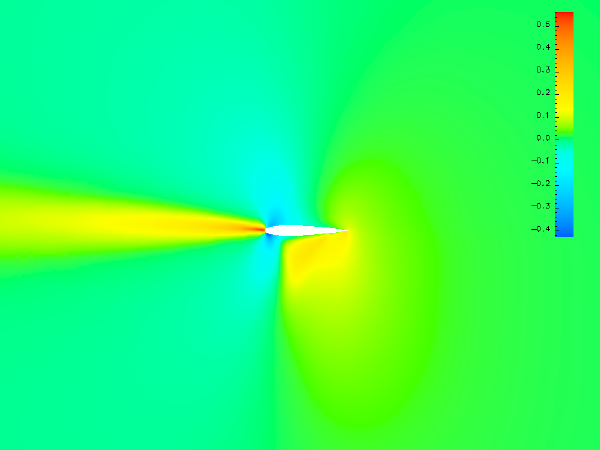}}
        \frame{\includegraphics[width=0.48\textwidth]{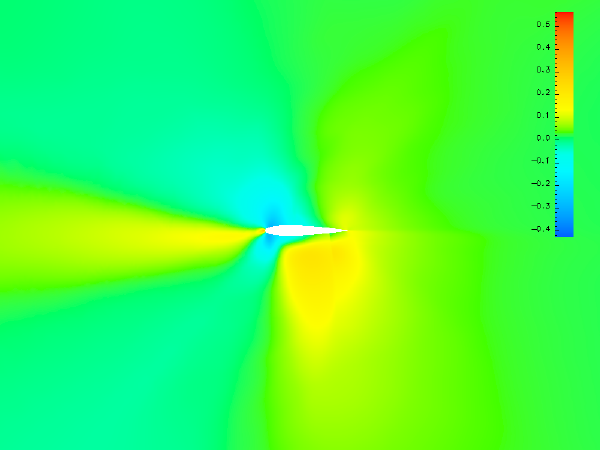}}   \\
        \frame{\includegraphics[width=0.48\textwidth]{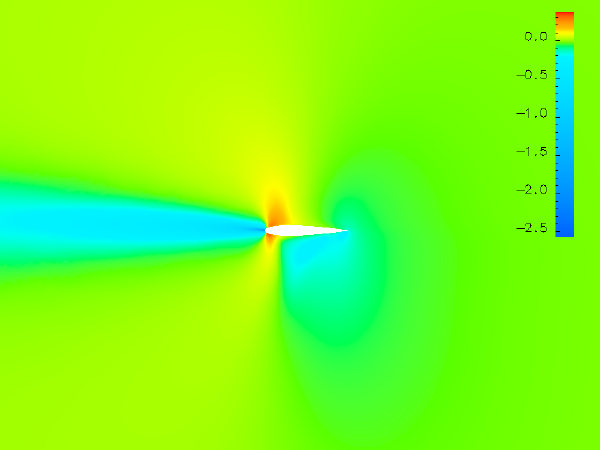}}
        \frame{\includegraphics[width=0.48\textwidth]{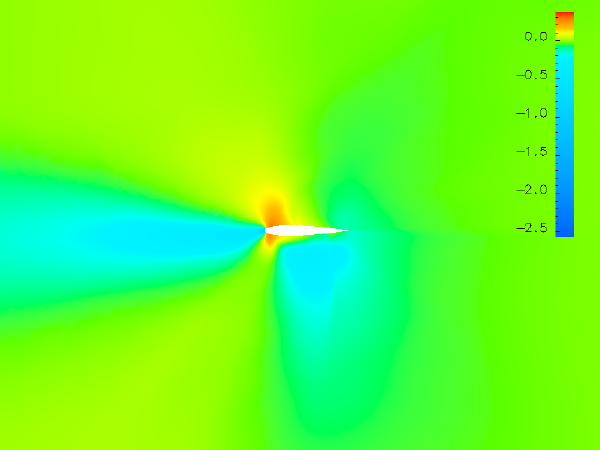}}\\
        \frame{\includegraphics[width=0.48\textwidth]{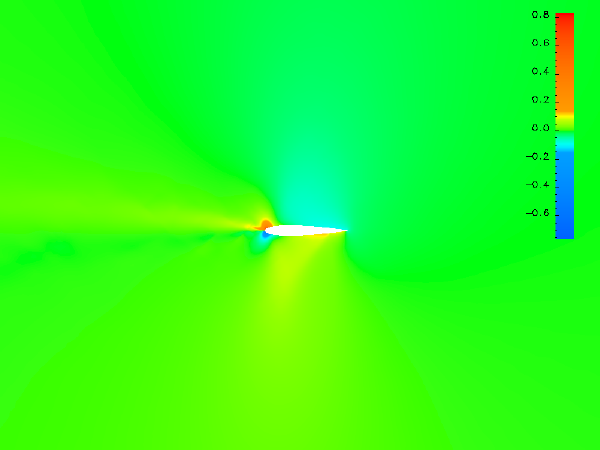}}
        \frame{\includegraphics[width=0.48\textwidth]{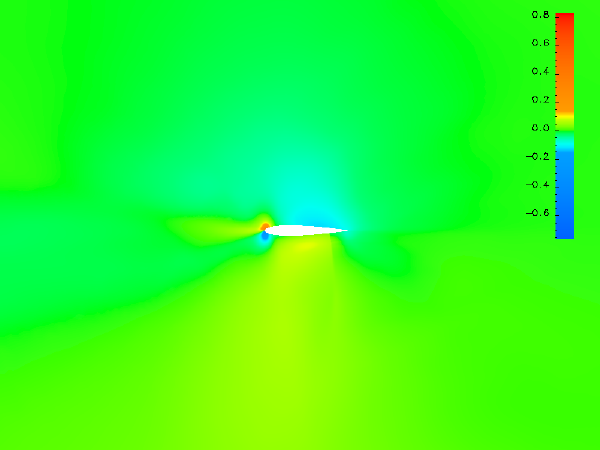}} \\  
        \frame{\includegraphics[width=0.48\textwidth]{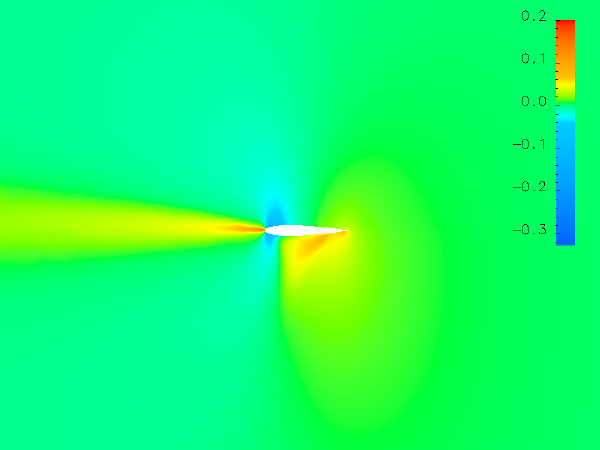}}
        \frame{\includegraphics[width=0.48\textwidth]{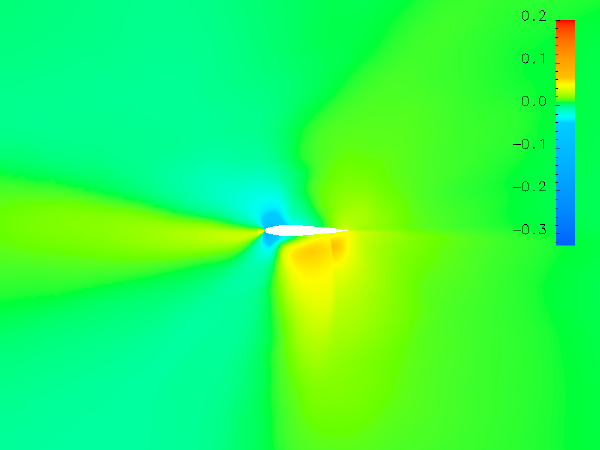}}
      \caption{Left column: The four dual variables generated by GMG solver; Right column: The four dual variables generated by CNNs model; NACA0012, Mach number 0.82, attack angle -1.65.}
      \label{dualCNN0.82solver}
      \end{figure}  

\begin{figure}[!h]\centering
  \includegraphics[width=0.45\textwidth]{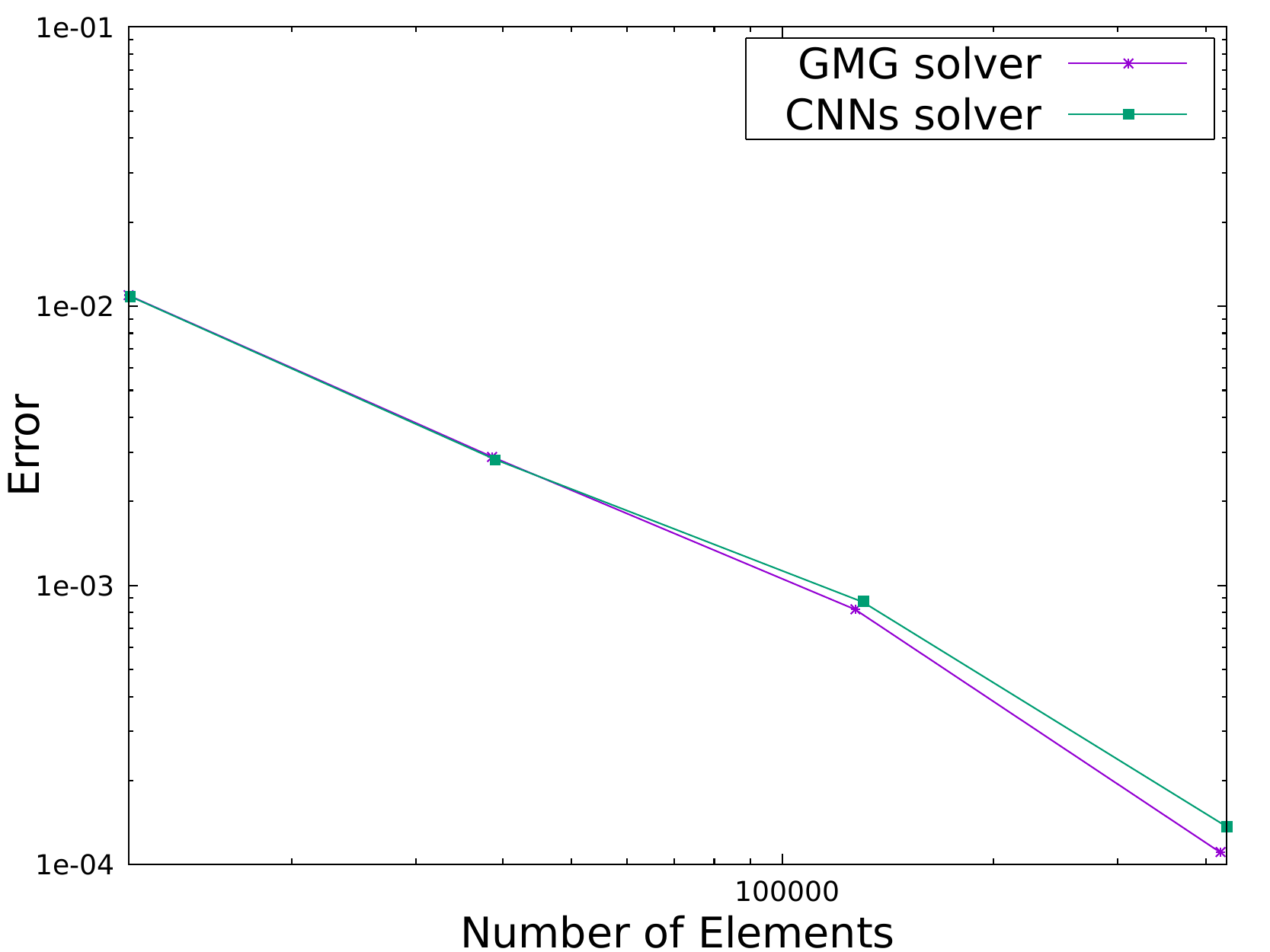}
  \includegraphics[width=0.45\textwidth]{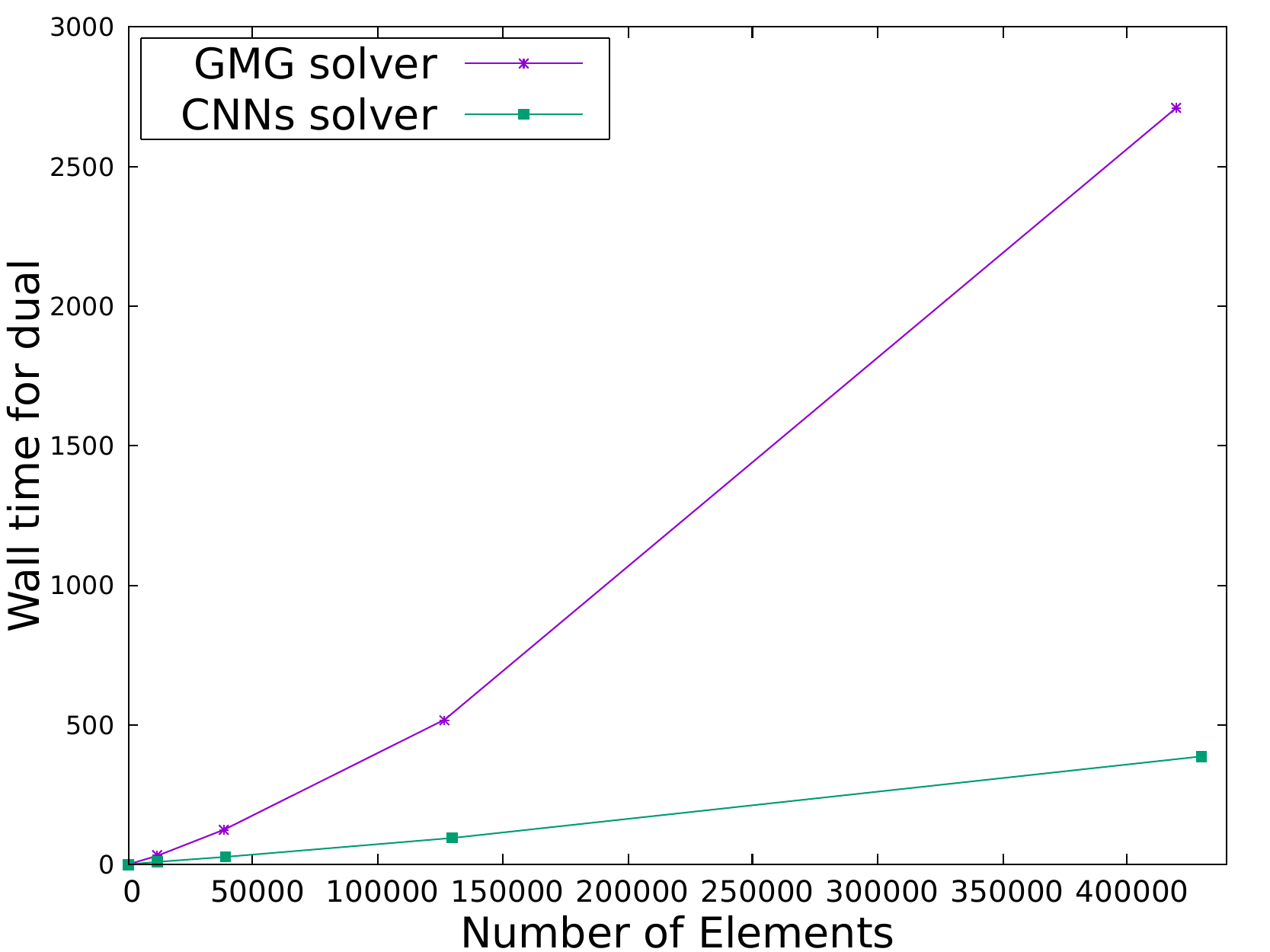}
\caption{Comparison of error and time cost with the dual solver of GMG form and CNNs form. NACA0012, Mach number 0.82, attack angle -1.65. }
\label{DUALvsGMG082}
\end{figure}

For fear of the convolutional neural networks may misunderstand the relation between primal and dual solver, we shall apply the model for other configurations. According to our experience, if the training sets are not sufficient enough, the CNNs dual solver may focus the refinement area on the top part of the airfoil. Then, we used the same ONNX model trained above to configurations with minus attack angle. Here we test the simulation with Mach number $0.82$, and attack angle $-1.65$. In Figure \ref{dualCNN0.82solver}, the trained model performs well on the different configurations, showing that the ability of our CNNs form dual solver can detect the main region for solving the target functional precisely. The precision and time spent for the adaptation are shown in Figure \ref{DUALvsGMG082}. The time spent with CNNs dual solver saves time for an order of magnitude while the precisions are comparable, which further verifies the reliability of our CNNs solver.

\begin{figure}[ht]\centering
  \includegraphics[width=0.45\textwidth]{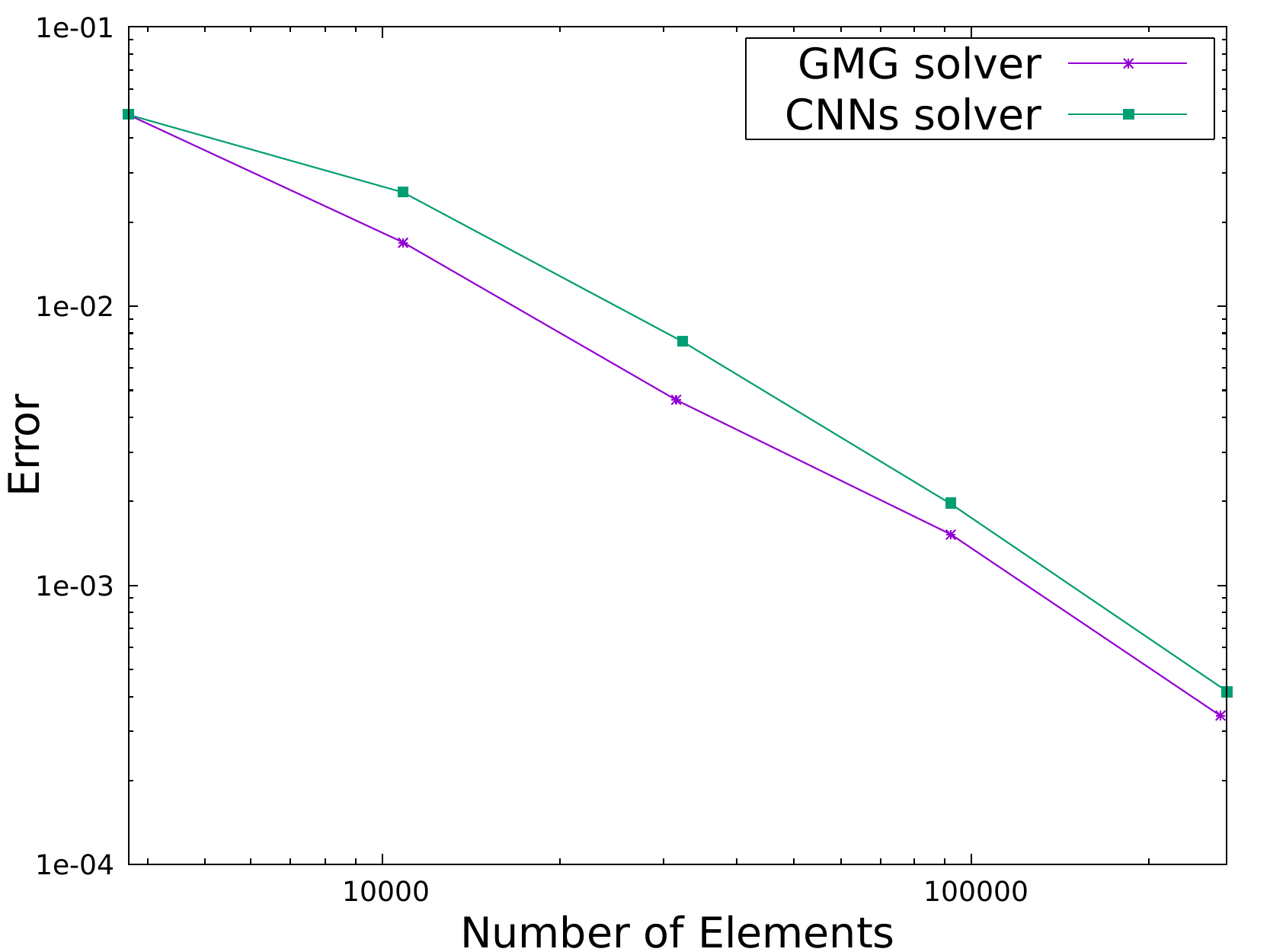}
  \includegraphics[width=0.45\textwidth]{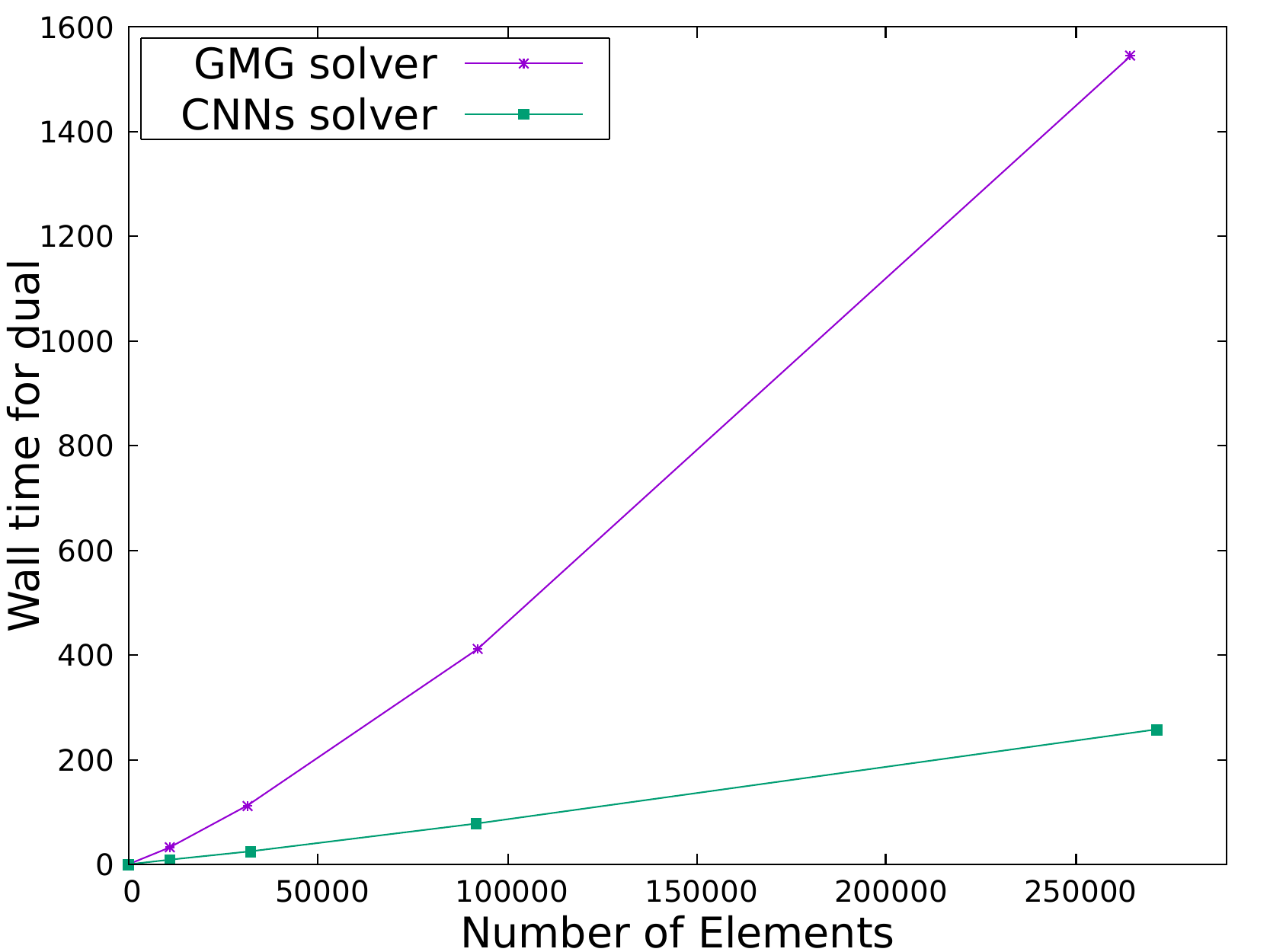}
\caption{Comparison of error and time cost with the dual solver of GMG form and CNNs form. RAE2822, Mach number 0.729, attack angle 2.31. }
\label{raecnn2822}
\end{figure}

\begin{figure}[!h]\centering
  \frame{\includegraphics[width=0.31\textwidth]{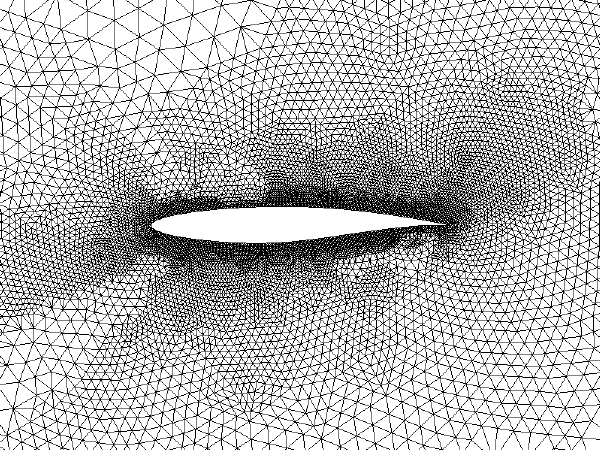}}   
  \frame{\includegraphics[width=0.31\textwidth]{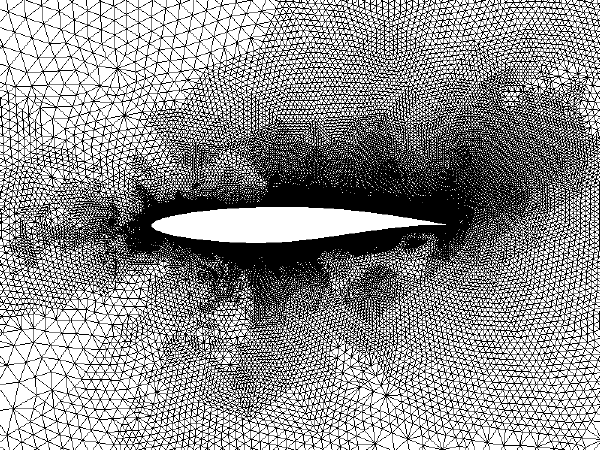}}
  \frame{\includegraphics[width=0.31\textwidth]{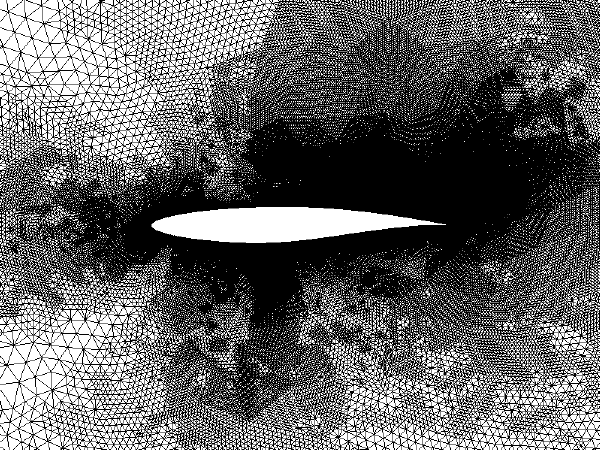}}
\caption{Consecutive refined meshes with CNNs dual solver. RAE2822, Mach number 0.729, attack angle 2.31$^\circ$.  }
\label{CNNadaptrae}
\end{figure}

Moreover, the trained model is tested on the RAE2822 airfoil. It demonstrates the model also performs satisfactorily since the precision is comparable to the GMG solver while saving time significantly as shown in Figure \ref{raecnn2822}. The consecutive refined meshes with the CNNs dual solver is illustrated in Figure \ref{CNNadaptrae}. The dual solver can also detect the main regions that influence the target functional most.

\subsubsection{Generalization on the multiple airfoils}
The benefit of our dual solver is the saving on the time for adaptation while the cost for the training is acceptable. In the training process, the dual consistency makes the model more reliable. Moreover, since the attention mechanism is introduced in this work. Even if the simulation is conducted on the unstructured mesh, it can still work well for multi-airfoil problems. In this part, we will show the feasibility of handling complicated geometry.

\begin{figure}[p]\centering
        \frame{\includegraphics[width=0.48\textwidth]{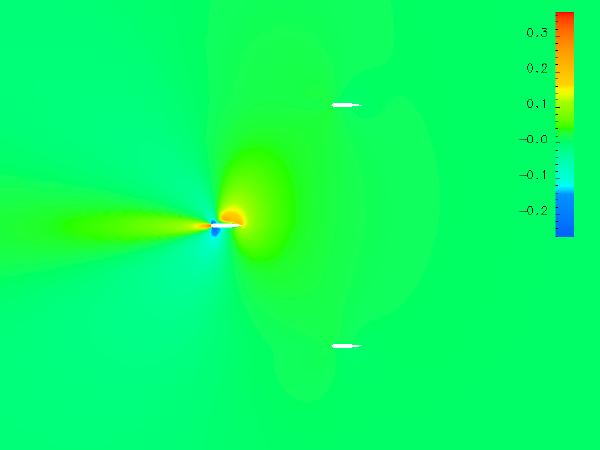}}
        \frame{\includegraphics[width=0.48\textwidth]{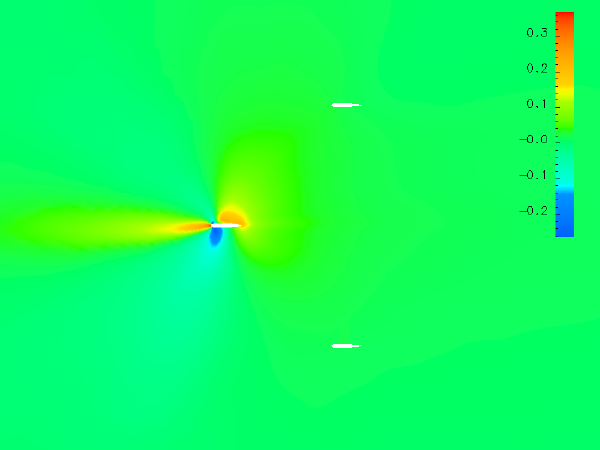}}\\ 
        \frame{\includegraphics[width=0.48\textwidth]{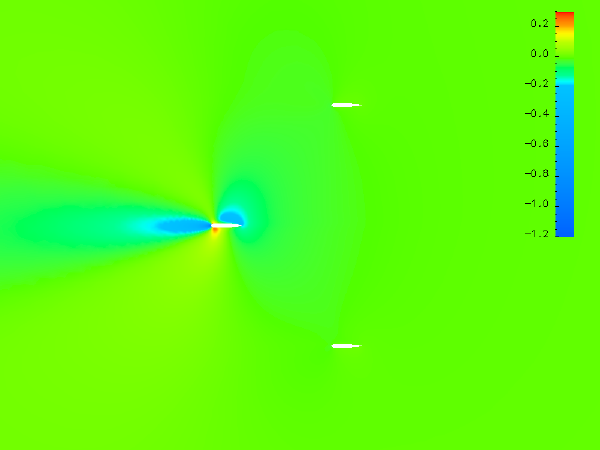}}
        \frame{\includegraphics[width=0.48\textwidth]{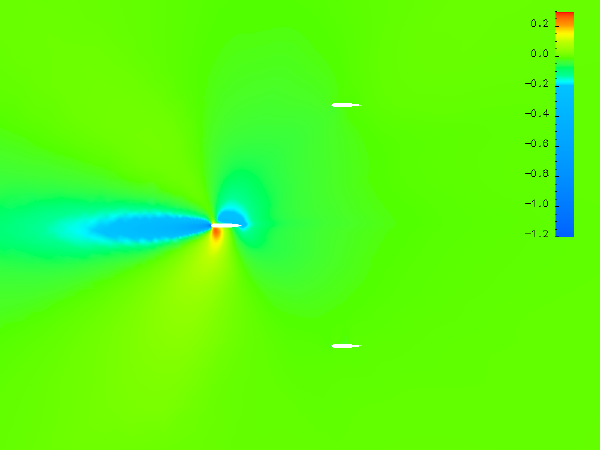}}\\
        \frame{\includegraphics[width=0.48\textwidth]{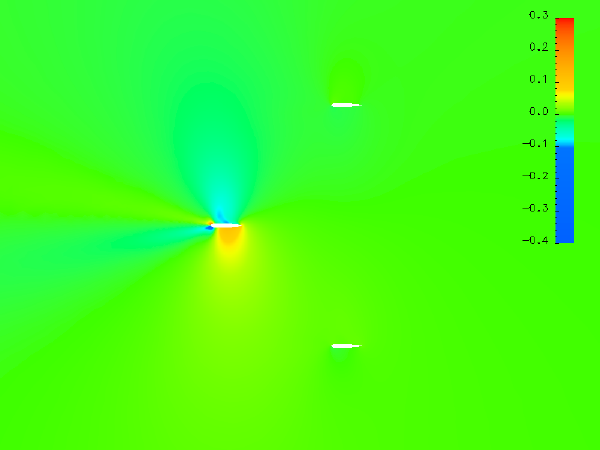}}
        \frame{\includegraphics[width=0.48\textwidth]{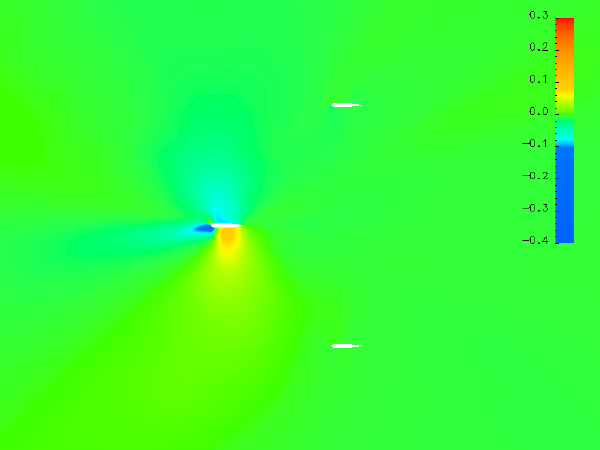}}\\ 
        \frame{\includegraphics[width=0.48\textwidth]{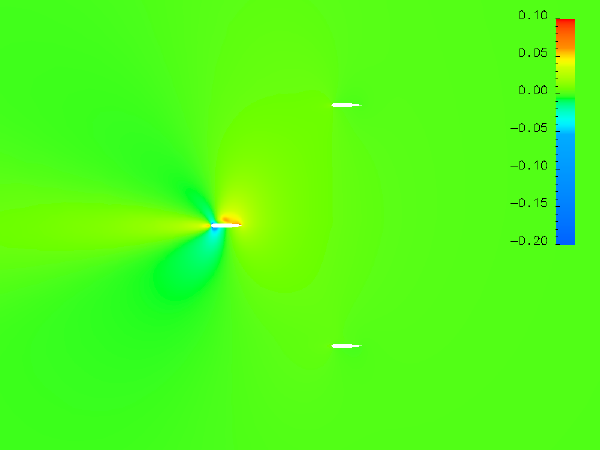}}
        \frame{\includegraphics[width=0.48\textwidth]{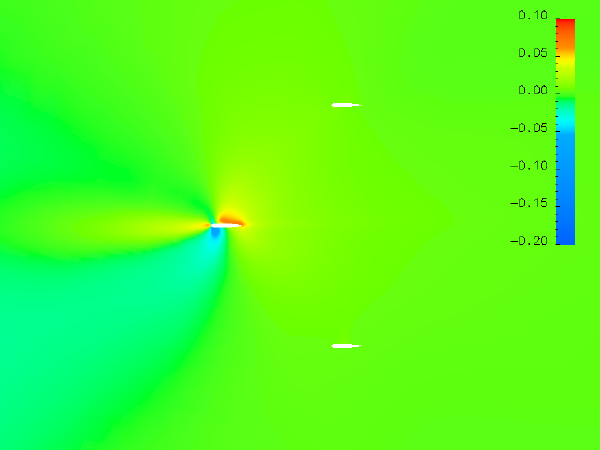}}
      \caption{Three-Airfoil model. Left column: The four dual variables generated by GMG solver; Right column: The four dual variables generated by CNNs model; Three NACA0012 airfoils, Mach number 0.8, attack angle 1.25$^\circ$.}
      \label{CNNThree}
      \end{figure}

We choose a model with 3 NACA0012 airfoils, and the target functional is the drag of the leading airfoil. According to the simulation in Figure \ref{CNNThree}, the CNNs form dual solver captured the main region and focused on the leading airfoil. 

\begin{figure}[!h]\centering
  \includegraphics[width=0.45\textwidth]{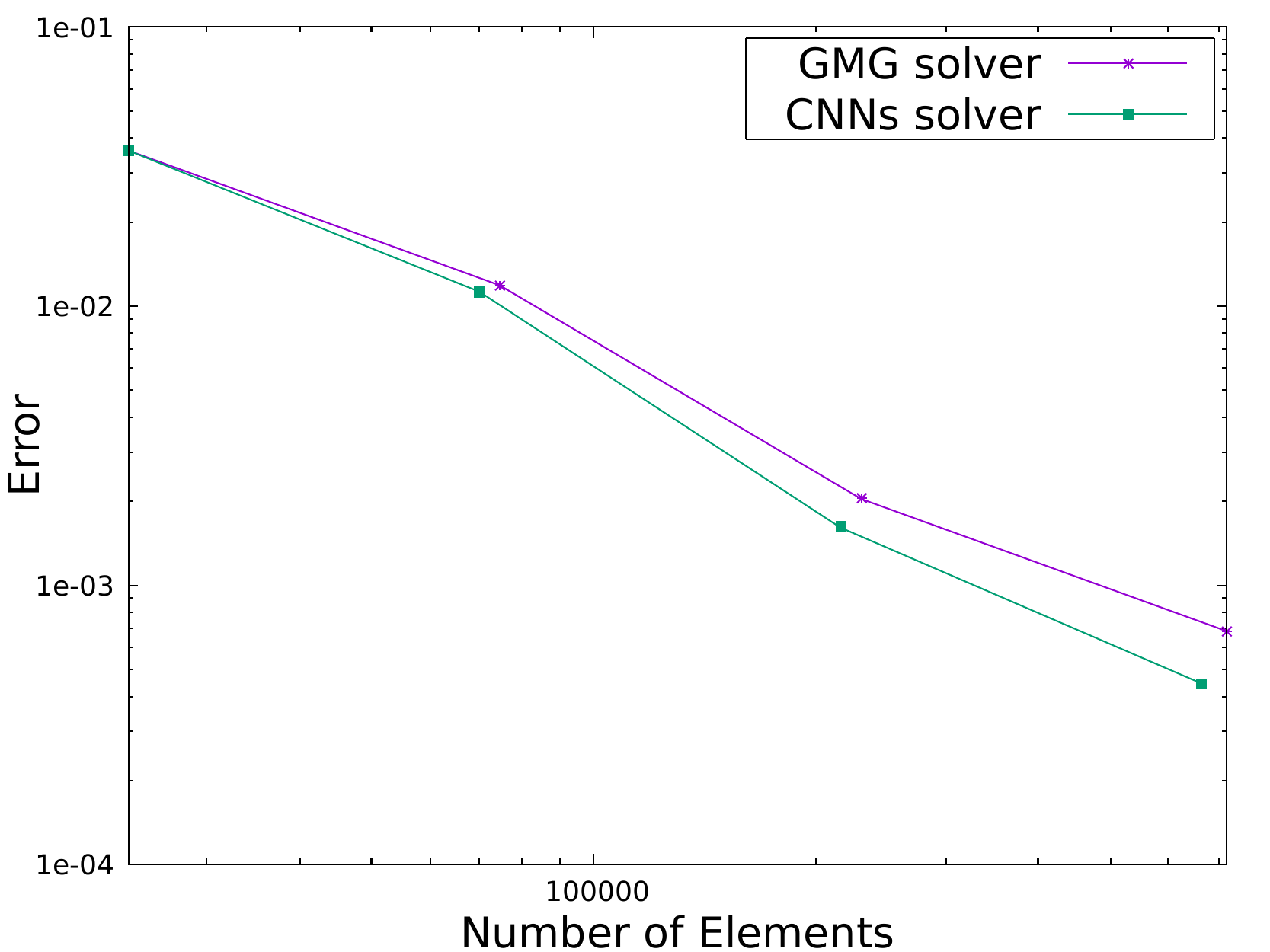}
  \includegraphics[width=0.45\textwidth]{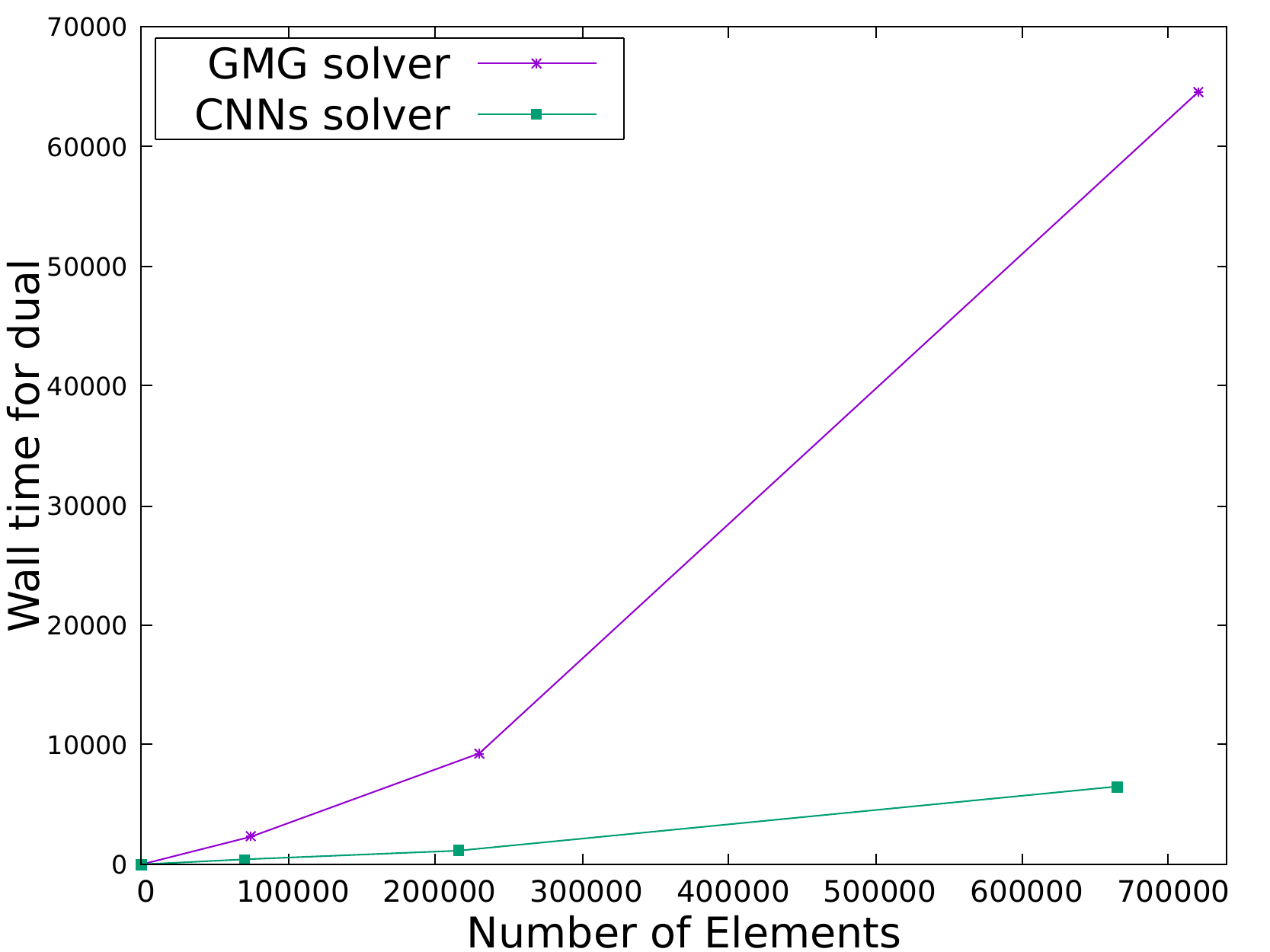}
\caption{Comparison of error and time cost with the dual solver of GMG form and CNNs form. Three NACA0012 airfoils, Mach number 0.8, attack angle 1.25$^\circ$. }
\label{DUALvsGMGThree}
\end{figure} 
\begin{figure}[!h]\centering
  \frame{\includegraphics[width=0.31\textwidth]{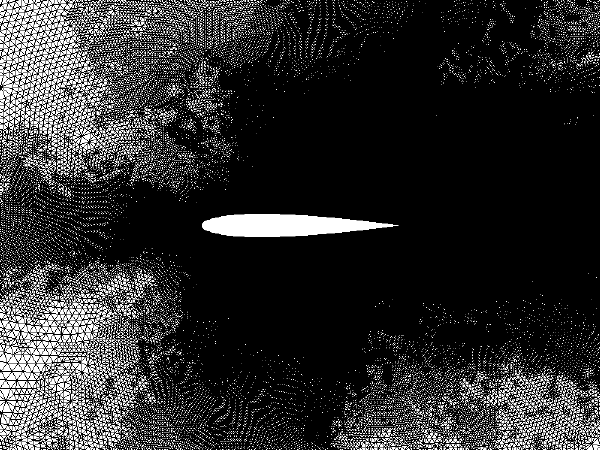}}   
  \frame{\includegraphics[width=0.31\textwidth]{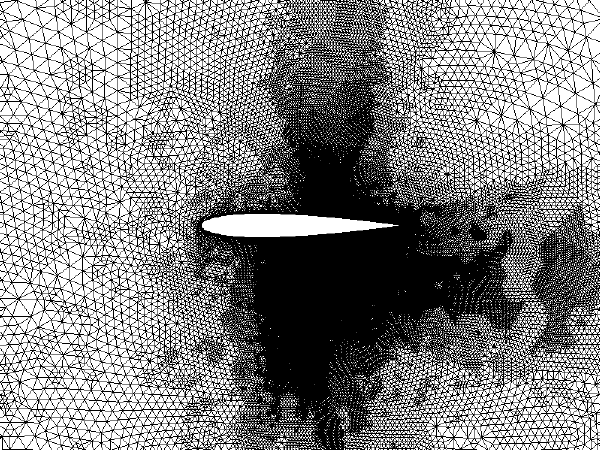}}
  \frame{\includegraphics[width=0.31\textwidth]{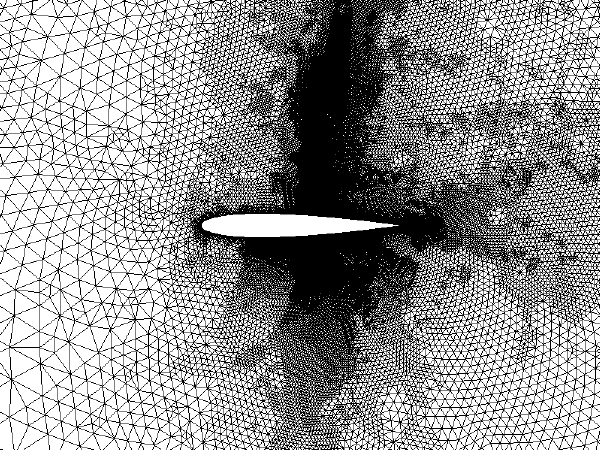}}
\caption{Meshes around the three different airfoils; Left: the leading airfoil; Middle: the upper airfoil; Right: the below airfoil.}
\label{CNNmeshThree}
\end{figure} 

It preserves a comparable precision similar to the GMG solver while saving time sharply as seen in Figure \ref{DUALvsGMGThree}. The initial mesh for this model is $23,466$. So the uniform refinement for this model takes a significant time for calculation. Then it is not suitable to compute the reference value with 5 times uniform refinement. For this sake, we only conduct the adaptation 3 times. As mentioned above, the training only needs data from the first two refinements. It will not cost too much computational resources. In this model, the precision from the CNNs form mesh generates a mesh whose quantity of interest is more accuracy than the GMG form solver.

The meshes around the different airfoils are shown in Figure \ref{CNNmeshThree}. The leading airfoil takes a balance between the residual and dual solutions. The refined areas not only centered around the leading edge but also around the upper half. However, the dual solutions are not distributed around the remaining two airfoils. Then the refined areas are only centered around the shock waves.

 \subsection{Software}
\subsubsection{AFVM4CFD}
AFVM4CFD is a library maintained by our group. It is an efficient solver which can solve the steady Euler equations with the $h$-adaptation method. Different modules such as k-exact reconstruction, parametric curves, and DWR-based $h$-adaptivity are integrated into this library. It is worth mentioning that with the AFVM4CFD, the Euler equations can be solved well
with a satisfactory residual that gets close to the machine precision. In this work, we further improved the efficiency by introducing the CNNs module. Now the shape optimal design is under construction. 

\subsubsection{Training Model}
The module build for the training can be seen from table \ref{tab:network_architecture}. The detail of the code can be found at https://github.com/ShanksFeng/CNNdualsolver.git. We will build more modules in the training in the future so that the performance will be better. The training can be conducted on different frameworks of the primal solver as long as the data structure is matched.

\section{Conclusion}
The dual-weighted residual-based mesh adaptation is very important for solving the quantity of interest in many practical issues. However, generating a suitable mesh for solving the target function is time-consuming. In this study, we constructed CNNs to predict the distribution of dual variables. It saves time in an order of magnitude compared with the traditional solver. More specifically, with time complexity of $O(n)$. Besides, to acclerate the whole algorithm, two different strategies are adopted. Firstly, a dynamic tolerance strategy is constructed in this work. It accelerates the calculation of the target functional and benefits for automatic selection without manual intervention. With such a strategy, a stable growth of mesh size in the adaptation process can be expected. Secondly, parallel computing technique is adopted for both the primal and dual solver, which enhance the computational efficiency markedly.
Moreover, the importance of dual consistency for generating reliable training datasets is validated in this work. The trained model has generalization capabilities that configurations beyond the trained categories can still be simulated precisely. It is worth mentioning that our convolutional neural networks can produce quality dual solutions on unstructured meshes.

Nonetheless, as the geometric information is not effectively integrated within the neural networks, unexpected oscillations may still occur in the vicinity of the boundary. This highlights the potential for future investigations utilizing physics-informed neural networks, aiming for a more accurate prediction of the dual variables. Furthermore, there is room for optimizing the architecture of our neural networks, with the potential to achieve more efficient adaptation in subsequent studies.
\bibliographystyle{unsrt} \bibliography{dualCNN}
\end{document}